\documentclass[12pt]{article}

\begin{document}

\input epsf.sty

\centerline{\bf S.P.Novikov\footnote{Sergey P. Novikov, IPST and
MATH Department, University of Maryland, College Park MD, USA and
Landau Institute, Moscow, e-mail novikov@ipst.umd.edu; This work
is partially supported by the Russian Grant in the Nonlinear
Dynamics. It is dedicated to the 70th birthday of Ya.G.Sinai. }}

 \vspace{0.3cm}

\centerline{\Large  Topology of the Generic Hamiltonian Foliations}
 \centerline{\Large on the Riemann Surfaces}

\vspace{0.5cm}

{\it Abstract.
  Topology of  the Generic Hamiltonian Dynamical Systems on the Riemann
   Surfaces given by the real part of the
generic  holomorphic 1-forms, is studied. Our approach is based on
the notion of Transversal Canonical Basis of  Cycles (TCB).
 This  approach
allows us to present a convenient combinatorial model of the whole
topology of the flow, especially effective for g=2. A maximal
abelian covering over the Riemann Surface
is needed here.
  The complete combinatorial model of the flow is
constructed. It  consists of the Plane Diagram and g straight line
flows in the  2-tori ''with obstacles''. The Fundamental Semigroup of positive
closed paths trasversal to foliation and topological invarians of trajectories are studied.
This work contains  an improved exposition of the results presented in the
authors recent preprint  \cite{NOV} and new results  calculating all TCB
in the 2-torus with obstacle, in terms of Continued Fractions.}
\vspace{0.5cm}

{\bf Introduction}. The family of parallel straight lines in the
Euclidean Plane $R^2$ gives  after factorization by the lattice
$Z^2\subset R^2$ the standard straight line flow in the 2-torus
$T^2$. It is a simplest ergodic system for the  irrational
direction. This system is Hamiltonian with multivalued Hamiltonian
function $H$ and standard canonically adjoint euclidean
coordinates $x,y$ (i.e. the  1-form $dH$ is closed but not exact,
and $x_t=H_y, y_t=-H_x$ where $H_x,H_y$ are constant). Every
smooth Hamiltonian system on the 2-torus without critical points
with irrational ''rotation number'' is diffeomorphic to the
straight line flow. Every $C^2$-smooth dynamical system on the
2-torus without critical points and with irrational rotation
number is $C^0$-homeomorphic to the straight line flow according
to the famous classical theorem. Various properties of the generic Hamiltonian Systems
on the 2-torus were studied in \cite{A}, ergodic properties were found in \cite{KhS}.\\

{\bf Question: What is going on
in the Riemann Surfaces with $g>1$?}\\
We consider  Hamiltonian Systems with multivalued
Hamiltonian $H$ (i.e.  given by the nonexact form $dH=0$ ).
 In
many cases we ignore time dependence of the trajectories and
discuss only the properties of  foliation $dH=0$ given by the
closed 1-form. {\it The present author started to investigate such
foliations  in early 1980s as a part of the newborn Topology of
The Closed 1-Forms. An important example was found in the Quantum
Solid State Physics describing the  motion of  semiclassical
electrons along the so-called Fermi Surface $M^2\subset  T^3$
for the single crystal
normal metals and low temperature in the strong magnetic field
(see \cite{N82}). The electron
trajectories exactly coincide with connectivity components of the
sections of Fermi Surface by the planes orthogonal to magnetic
field. One might say that they are the levels of quasiperiodic
function on the plane with 3 periods. An extensive study of that
class was performed by the author's Moscow Topological Seminar
since early 1980s and was continued in Maryland.
Remarkable ''topological complete
integrability''  and ''topological resonance'' properties in the
nonstandard sense were discovered here, for the set of directions of magnetic field of
the full measure on the 2-sphere $S^2$. These properties play a
key role in the physical applications \cite{DN05}.}

 Following  class of generic Hamiltonian systems on the Riemann Surfaces
 was studied  by many
 mathematicians (no applications outside of pure
mathematics were found for them  until now, unfortunately): Take
any nonsingular compact Riemann Surface $V$, with genus  $g$.
  Every holomorphic 1-form
$\omega\in C^g$ defines a Hamiltonian system (foliation) $\Re$ on
the manifold $V$: $$\Re=\{ \omega^R=0\}$$ because $d(\omega^R)=0$.
Here $\omega=\omega^R+i\omega^I$. Most people considered a more general class
of ''foliations $\Re$ with invariant transversal measure'' on
Riemann Surfaces. Every holomorphic quadratic differential $\Omega$
 defines such foliation $\Re$ by the formula $(\sqrt{\Omega})^R=0$.
This is a locally hamiltonian foliation  (i.e. it admits a
transversal measure) but  non-orientable (it does not admit global
 time direction  if $\Omega\neq \omega^2$). We do not consider
 nonorientable case. It does not define  Hamiltonian Dynamical System
 (even locally it has different generic singularities).

The  systems with transversal invariant measure were studied since
early 1960s by the following method (see in the book \cite{KHa}):
Take any closed curve $\gamma$ transversal to our foliation $\Re$.
Assume that almost every nonsingular trajectory is dense. For
every point $Q\in \gamma$ except finite number
the trajectory started in the point $Q$
returns to the curve $\gamma$ first time at the new point $P$ (in positive
direction of time). We define a ''Poincare map'' $Q\rightarrow P$.
This map  preserves  transversal measure equal to the
restriction of the form $dH$ on the curve $\gamma$.   {\bf It is
the Energy Conservation Law} for  Hamiltonian System. So our
transversal closed curve $\gamma$ is divided into
 $k=k_{\gamma}$ intervals $\gamma=I_1+I_2+\ldots+I_k$. The
Poincare map looks  as a  permutation of $k$ intervals
on the circle; this map is ill-defined in the finite number of points
only. The time $t_Q$ varies continuously in the interior of each interval
$Q\in I_j$.

No doubt, the use of  closed transversal curves is extremely
productive. At the same time,  we are not satisfied by this
approach; Following questions can be naturally asked:

1.This method essentially ignores time/length and homology/homotopy
classes of trajectories
starting and ending in $\gamma$. Our goal is to present some sort of global topological
description  of the flow (or foliation $\Re$) on the algebraic
curve $V$   similar to the case of
genus 1 as much as possible.

 2.There are many different closed transversal curves
in the foliation $\Re$. How this picture depends on the choice of
tranversal curve $\gamma$? Nobody classified them yet as far as I
know. Indeed, in the theory of codimension 1 foliations developed
by the present author in 1960s (see\cite{N65}) several algebraic
structures were defined for the closed transversal curves:

Consider all closed
positively (negatively) oriented transversal curves starting and
ending in the point $Q\in V$. We can multiply them. Transversal homotopy classes
of such curves generate {\bf A Transversal Fundamental Semigroup}
$\pi_1^+(\Re,Q)$ and its natural homomorphism into the fundamental
group (even, into the fundamental group of the unit tangent
$S^1$-bundle $L(V)$
$$\psi^{\pm}:\pi_1^{\pm}(\Re,Q)\rightarrow\pi_1(L(V),Q)\rightarrow\pi_1(V,Q)$$
The set of all closed transversal curves naturally maps into the
set of conjugacy classes in fundamental group. We denote it also
by $\psi$. The Transversal Semigroups might depend of the leaf
where the initial point $Q$ is chosen.

{\bf How to calculate these invariants for the Hamiltonian Systems
 on the algebraic curves?}

The simplest fundamental properties of these foliations are following:

Property 1. They have only saddle type critical points (no centers). In the
generic case such foliation has exactly 2g-2 nondegenerate
saddles.

Property 2. Every nonempty closed  transversal curve $\gamma$ is
non-homologous to zero. Every period of the form $dH$ is positive
$\oint_{\gamma}dH>0$ for the positive transversal
curve. So the composition $$\pi_1^+(\Re,Q)\rightarrow
\pi_1(V,Q)\rightarrow H_1(V,Z)\rightarrow R$$ does not map any
element into zero. Every periodic trajectory $\gamma$ is such that
$\oint_{\gamma}dH=0$. Therefore for the generic case every periodic
trajectory is homologous to zero. We call systems without periodic trajectories
{\bf Irreducible}.

{\bf Definitions}

a.We say that foliation on the Riemann surface $V$   belongs to
the  class $T$  if there exist a {\bf Transversal Canonical
Basis} of curves $$a_1,b_1,\ldots,a_g,b_g$$ such that
all these curves are non-selfintersecting (simple), transversal
to  foliation, and the curves $a_j$ and $b_j$ transversally cross each
other
 exactly in one point. Other curves do not cross
each other.

b.We say that foliation  belongs to the class $T^{0}$  if
there exists a canonical basis such that all $a$-cycles
$a_1,\ldots,a_g$ are simple, do not cross each other
and are transversal to  foliation.
We say that foliation  belongs to the {\bf Mixed Class}  of the
Type $T^{k},0,k=1,\ldots,g$,  if there exists  an incomplete canonical basis
$a_j,b_q,j\leq g,q\leq k,$ such that all these cycles are transversal to foliation,
simple, and  only pairs crossing each other are $a_j$
and $b_j$ for $j=1,2,\ldots,k$. All intersections  are transversal
 and consist of one point each. The complete case is
 $T=T^{g}$

\newtheorem{rem}{Remark}

\begin{rem}
A number of people including Katok, Hasselblatt, Hubbard, Mazur,
Veech,
 Zorich, Konzevich, McMullen  and others wrote a lot of works related to study the
ergodic properties  of foliations with  ''transversal measure'' on
the Riemann Surfaces, and their total moduli space
 (see  \cite{K,Sat,KHa,HM,KMS,M,V,MT,Z,Z1,Z2,ZK,Kn,F,McM}).  People investigated
recently  closed geodesics of the flat Riemannian Metric
$ds^2=\omega\bar{\omega}$ singular in the critical points of the
holomorphic 1-form $\omega$. These geodesics consists of all
trajectories of the  systems
$(\exp\{i\phi\}\omega)^R=0$ on  $V$.  Closed
geodesics appear for the special  number of angles $\phi_n$ only.
Therefore we don't have them describing the generic Hamiltonian
systems. Beautiful analogs of the Poincare 4-gons
associated with  flat metrics
 with singularities were invented  and used.
 Our intention is to describe
topology of  the generic foliations
given by the equation $\omega^R=0$ for the holomorphic form
$\omega$. As A.Zorich pointed out to us, people
 certainly observed  some   features which
may illustrate our ideas (see
 \cite{KHa,McM}).
 Our key idea of transversal canonical
 basis did not  appeared before.
\end{rem}

\vspace{0.3cm}

 \centerline{\Large Section 1. }

\centerline{\Large  Examples: The Real Hyperelliptic Curves.}

\vspace{0.5cm}

Let us show the examples of
foliations in the classes $T,T^{k}$. Consider  any real
hyperelliptic curve of the form
$$w^2=R_{2g+2}(z)=\prod_j (z-z_j),z_j\neq z_l$$ where all roots
$z_j$ are real and ordered naturally $z_1<z_2<\ldots <z_{2g+2}$.

\newtheorem{lem}{Lemma}
\begin{lem}
For every complex non-real and non-imaginary number $u+iv,u\neq
0,v\neq 0$, generic polynomial $P_{g-1}(z)$ and  number $\epsilon$
small enough, the Hamiltonian system defined by the closed
harmonic 1-form below belongs to the class  $T^{0}$. All cycles
$a_1,\ldots,a_g,c_1,\ldots,c_{g+1}$ are transversal to foliation
and cross each other transversally,
$$ \omega=\frac{(u+iv+\epsilon P_{g-1}(z))dz}{\sqrt{R_{2g+2}(z)}},
\omega^R=0$$

 For $g=1,2$ this foliation belongs to the class $T$.
For $g>2$ it belongs to the class $T^{2}$.

\end{lem}
Proof. The polynomial $R=R_{2g+2}(z)$ is real on the cycles
$a_j=p^{-1}[z_{2j}z_{2j+1}],j=1,...,g$ and purely imaginary on the
cycles $c_q$ located on the real line $x$ immediately before and
after them $c_q=p^{-1}[z_{2q-1}z_{2q}],q=1,2,\ldots,g+1$.

 For $g=1$ we take a canonical basis
$a_1,b_1=c_1$. For $g=2$ we take a canonical basis
$a_1,a_2,b_1=c_1,b_2=c_3$. We shall see that all these cycles are
transversal to foliation.

  For
$\epsilon=0$ we have $$\omega^R=(udx-vdy)/R^{1/2}=0,y=0,x\in a_j$$
$$\omega^R=(vdx+udy)/iR^{1/2}=0, y=0,x\in c_j$$ In both cases we
have for the second component $dy\neq 0$ for the direction of
Hamiltonian system. So the transversality holds for all points on
the cycles except (maybe) of the branching points $z_j$. The
cycles $c_l,a_k$ are orthogonal to each other in the crossing
points (i.e. in the branching points).  We need to check now that
in all branching points $z_j$ the angle between trajectory and
both cycles $a,c$ is never equal to $\pm\pi/2$ because they meet
each other with this angle exactly. After substitution
$w'^2=z-z_j$ we have
$$\omega=(u+iv)2w'dw'/w'F_j^{1/2}=2(u+iv)dw'/F^{1/2}_j$$ where
$F_j=\prod_{l\neq j} (z-z_l)$. The real function
$F_j(x)$ for real $x\in R$  does dot
change sign passing through the point $x=z_j\in R$:  $F_j(z_j)\neq 0$.
For $w'=f+ig$ we have $\omega^R=2(udf-vdg)/F_j^{1/2}$ if
$F_j(z_j)>0$ or $\omega^R=2i(udg+vdf)/F^{1/2}$ if $F_j(z_j)<0$. We
have $x=f^2-g^2,y=2fg$. The condition $y=0$ implies the equation
$fg=0$, i.e. the union of the equations $f=0$ and $g=0$. It is
exactly an orthogonal crossing of our cycles. In both cases our
system $\omega^R=0$ implies transversality of trajectories to
both cycles $df=0$ and $dg=0$ locally.

We can choose now the $a$-cycles as any subset of $g$ cycles out of $a,c$
not crossing each other, for example $a_1,\ldots,a_g$. All of
them are transversal to foliation. Therefore we are coming to the
class $T^0$. Now we choose following two $b$-cycles: $b_1=c_1$ and
$b_g=c_{g+1}$. According to our arguments, this choice leads to
the statement that our foliation belongs to the class $T$ for
$g=1,2$ and  to the class $T^2$ for all $g\geq 2$.

 Lemma
is proved now because  small $\epsilon$-perturbation cannot
destroy transversality along  the finite family of compact cycles.

We choose now the generic perturbation such that {\bf all critical
points became nondegenerate, and all saddle connections and
periodic trajectories nonhomologous to zero disappear.}

Let now $R=R_{2g+2}=\prod_{j=1}^{2g+2}(z-z_j)$ is a polynomial of
even degree as above with real simple roots $z_j\in R$, and
$\omega=P_{g-1}(z)dz/R^{1/2}$ is a generic holomorphic 1-form. Let
$P=u+iv$ where $u,v$ are real polynomials in the variables $x,y$.
Only the zeroes $v=0$ in the segments
$[z_{2q-1}z_{2q}],q=1,...,g+1$, and the zeroes $u=0$ in the
segments $[z_{2j}z_{2j+1}],j=1,...,g$, are important now. We
assume that $u(z_k)
\neq 0$ and $v(k)\neq 0$ for all
$k=1,\ldots,2g+2$.
\begin{lem}  Remove all open segments  containing the important zeroes, i.e. the
open segments $[z_{2q-1}z_{2q}]$ containing the zeroes  $v=0$, and
all open segments $[z_{2j}z_{2j+1}]$ containing the zeroes $u=0$.
If remaining segments are enough for the construction of the
half-basis $a_1,\ldots,a_g$ (i.e. there exists at least $g$
disjoint closed segments between them), then our foliation belongs
to the class $T^0$. In particular, it is always true for $g=2$
where $u$ and $v$ are the linear functions (and have no more than
one real zero each).
 Let there is no real  zeroes $v=0$
  in the segments $[z_1z_2],[z_5z_6]$,
and no real zeroes $u=0$ in the segments $[z_2z_3],[z_4z_5]$  for
$g=2$; Then this foliation belongs to the class $T$.
 Let $g>2$, all real   zeroes $v=0$ belong to the open intervals
$$(-\infty,z_1),(z_3z_4),(z_5z_6),\ldots,
 (z_{2g-1}z_{2g}),( z_{2g+2},+\infty)$$ and all real zeroes $u=0$
are located   in the open segments
$$(-\infty,z_2),(z_3z_4),...,(z_{2g-1}z_{2g},(z_{2g+2},+\infty)$$
In this case  foliation $\omega^R=0$ belongs to the class $T^2$.

\end{lem}
Proof is exactly the same as above.

For every foliation of the class $T^0$ we cut Riemann Surface $V$
along the transversal curves $a_j$. The remaining manifold
$\tilde{V}$ has a boundary
$\partial{\tilde{V}}=\bigcup_ja_j^{\pm}$ where $a_j^{\pm}=S^1$;
the foliation enters it from inside for $(a^-_j)$ and
 from outside for $(a_j^+)$.
The domain $\tilde{V}=D^2_*\subset R^2$ is a disc with
 $2g-1$ holes removed from  inside: the external boundary
is exactly $a_1^+$. The boundaries of inner holes
are $a_1^-,a_j^{\pm}, j>1$.
All boundaries are transversal to our system (see Fig 1) with periods
 $$\oint_{a_j}\omega^R=|a_j|>0, |a_1|\geq
|a_2|\geq \ldots \geq |a_g|>0$$
 The system does not have critical
points except $2g-2$ nondegenerate saddles. Every trajectory
starts at the {\bf in-boundary} $\bigcup a^+_j$ and ends at the
{\bf out-boundary} $\bigcup a_j^-$  (see Fig 4).

Consider a maximal abelian $Z^{2g}$-covering  $V'\rightarrow V$
 with basic shifts $a_j,b_j:V'\rightarrow V'$
for every class $T$ foliation. Its fundamental domain can be
obtained cutting $V$ along the Transversal Canonical Basis. The
connected pieces $A_j$ of the boundary exactly represent free
abelian groups $Z^2_j$ generated by the shifts $a_j,b_j$. Every
 boundary component $A_j=A^+_j\bigcup A^-_j$  looks
like the standard boundary square  on the plane (see Fig 2). The foliation
near the boundary $A_j$ looks like  standard straight line flow
at the space $R^2_j$ with lattice $Z^2_j$. Our fundamental domain
$\bar{V}\subset V',\partial \bar{V}=\bigcup A_j$, is restricted to
the inner part of 2-parallelogram in every such plane $R^2_j$.
Topologically this domain  $\bar{V}$ is a 2-sphere with $g$ holes
(squares) , with
 boundaries $A_j$.

Algebraic Geometry  constructs $V'$ analytically  (the Abel Map):
$$A=(A^1,\ldots,A^g):V\rightarrow V'\subset C^g,
A^j(Q)=\int_{P}^Q\omega_j,j=1,\ldots,g$$ Here
$\omega_k$ is a normalized basis of
holomorphic forms:
$$\oint_{a_j}\omega_k=\delta_{jk},\oint_{b_j}\omega_k=b_{jk}=b_{kj}$$
The form $\omega=\sum u_k\omega_k$ is generic here. It defines a
one-valued function $F$:
$$F:V'\rightarrow
C;  F(Q')=A^{\omega}(Q')=\sum u_kA^k(Q')$$ where $Q'$ corresponds to
$Q$ under the projection $V'\rightarrow V$.  The basic shifts are
$a_j,b_j:V'\rightarrow V'$:
$$F(a_j(Q'))=F(Q')+u_j,F(b_j(Q'))=F(Q')+\sum_k u_kb_{kj}$$ The
levels $F^R=const$ are exactly the leaves of our foliation
$\omega^R=0$ on the covering $V'$. The map $\pi_1(V)\rightarrow
H_1(V,Z)=Z^{2g}\rightarrow R$ is defined by the correspondence:
$$a_j\rightarrow u_k^R,b_j\rightarrow (\sum_ku_kb_{kj})^R$$ The
map $\psi: \pi^+_1(\Re)\rightarrow \pi_1(V)\rightarrow R^+$ of the
positive transversal semigroup (above) certainly belongs to the
$(n,m)\in Z_+^{2g}(\omega^R)$ if
$$\omega^R(\sum_{m_k,n_l}m_ka_k+n_lb_l)>0$$

\vspace{0.3cm}

\centerline{\Large Section 2. Some General Statements.}

\vspace{0.3cm}

Let us introduce  the class of {\bf  Positive Almost Transversal Curves}. They are
piecewise smooth and consist of  the smooth pieces
 of the First and Second Type: The
First Type pieces are positive and transversal to foliation.
Every Second Type piece belongs to one trajectory.

 {\bf Every almost
transversal curve can be approximated by the smooth transversal
curve} with the same endpoints (if there are any). In many cases
below we construct closed almost transversal curves and say
without further comments that we constructed smooth transversal
curve.

The {\bf Plane Diagram} of foliation of the type
$T^k$ with incomplete transversal canonical basis is  {\bf Topological Type}
of foliation on the plane domain $\bar{V}$ obtained from $V$ by the
cuts along this basis (see Fig 10).

\newtheorem{thm}{Theorem}

\begin{thm} Every generic foliation given by the holomorphic 1-form $\omega^R=0$
on the algebraic curve of genus 2
admits  Transversal Canonical Basis
\end{thm}

\begin{lem}Every generic foliation $\omega^R=0$ belongs to the
 class $T^0$ for genus equal to  2
\end{lem}
Proof. Take any trajectory such that its limiting set in both
directions  contains at least one nonsingular point.
 Nearby of the limiting
nonsingular point our trajectory appears infinite number of times.
Take two nearest returns and join them by the small
transversal segment. This closed curve consists of
the piece of trajectory  and small transversal segment. It is an
almost transversal curve. Therefore it can be approximated by the closed
non-selfintersecting smooth transversal curve. Take this curve
as a cycle $a_1$. Now  cut $V$ along this curve and get the
surface $\tilde{V}$ with 2 boundaries
$\partial\tilde{V}=a_1^+\bigcup a_1^-$. Take any trajectory
started at the cycle $a_1^+$ and ended at $a^-_1$. Such trajectory
certainly exists. Join the ends  of this trajectory on the cycle
$a^-_1$ by the positive
 transversal segment  along the cycle $a_1$ in $V$. We
get a transversal non-selfintersecting cycle $b_1$, crossing $a_1$
transversally in one point. Cut  $V$ along the pair $a_1,b_1$.
We get a square $\partial D^2=a_1b^{-1}_1a_1^{-1}b_1\subset R^2$
with   1-handle attached to the disk $D^2$ inside. Choose
such notations for cycles  that our foliation enters this
square along the piece $A_1^+=a_1b^{-1}_1\subset
\partial D^2$, and leaves it along the piece
$A_1^-=a_1^{-1}b_1$ (see Fig 3). Nearby of the angle  these
pieces $R^{\pm}$ are attached to each other. The trajectories here
spend small time inside of the square: they enter and leave it cutting the angle
(see Fig 3). Move  inside from the  ends of the segment
$A_1^+$. Finally we find 2 points $x^+_{1,1},x^+_{2,}\in A_1^+$ where
this picture ends (because topology is nontrivial inside of the square).
 These points
are the ends of the separatrices of two different saddles
 $x^+_{1,1}\neq x^+_{2,1}$.  Take
any point $y$ between $x^+_{1,1}$ and $x^+_{2,}$ (nearby of the left one
$x^+_{1,1}$). The trajectory passing through $y$ crosses $A_1^-$
in the point $y'$ (see Fig 3). Join  $y'$ with point $\tilde{y}$ equivalent to $y$
 by the transversal segment $\sigma$ along the curve
$A_1^-$ in positive direction.
 This is a closed curve $c_1$ almost transversal to our
foliation. If initial point is located on the cycle $a_1$, the
curve $c_1$ does not cross the cycle $b_1$. We take cycles
$b_1,c_1$ as a basis of the transversal $a$-cycles. If initial
point is located on the cycle $b_1$, the curve $c_1$ does not
cross $a_1$. In this case we take $a_1,c_1$ as a basis of the
transversal $a$-cycles. The cycle $c_1$ is
non-homologous to $b_1$ in the first case. Therefore it is a right
basis of the $a$-cycles which is transversal. The second case is
completely analogous. Our lemma is proved.

\begin{lem} For every generic foliation $\omega^R=0$ of the class $T^0$
on the algebraic curve of genus 2,
 the transversal basis of $a$-cycles can be extended to the full
Transversal Canonical Basis $a,b$, so every $T^0$-class foliation
belongs to the class $T$.
\end{lem}

Proof. Cut the Riemann Surface $V$ along the transversal cycles
$a_1,a_2$. Assuming that $\oint_{a_1}\omega^R=|a_1|\geq
|a_2|=\oint_{a_2}\omega^R$, we realize this domain as a plane
domain $D^2_*$ as above (see Fig 4 and the previous section). Here
an external boundary $\partial_{ext} D^2_*=a_1^+$ is taken as our
maximal cycle. The elementary qualitative intuition shows that
there are only two different topological types of the plane
diagrams (see Fig 4, a and 4,b): The {\bf first case} is
characterized by the property that for each saddle all its
separatrices end up in the four different components of boundary
$a_1^{\pm},a_2^{\pm}$. We have following matrix of the trajectory
connections of the type $(k,l):a^+_k\rightarrow a^-_l$  for the
in- and out-cycles and their transversal measures:

 $a^+_1\rightarrow a_1^-$ (with measure $a$), $a^+_k\rightarrow a^-_l, k\neq l$
(with measure $b$ for $(k,l)=(1,2),(l,k)=(2,1)$), and
$a^+_2\rightarrow a^-_2$ with measure $c$. All measures here are
positive. We have for the measures of cycles:
$|a_1|=a+b,|a_2|=c+b$. This topological type does not have any
degeneracy for $a=c$.

The diagonal trajectory connections of the type $(l,l)$ generate
the transversal $b$-cycles closing them by the transversal pieces
along the end-cycles in the positive direction.

In the  second case we have following matrix of trajectory
connections:

$a_1^+\rightarrow a^-_1$ with measure $a>0$, $a^+_1\rightarrow
a^+_2$ with measure $b>0$, $a_2^+\rightarrow a_1^-$ with the same
measure $b>0$. We have $b=|a_2|,a+b=|a_1|$. So we do have the
trajectory-connection $a_1^+\rightarrow a_1^-$, but we do not have
the second one, of the type $(2,2)$. However, we may connect
$a_2^+$ with $a_2^-$ by the almost transversal curve as it is
shown in the Fig 4,b), black line $\gamma$.  So we construct the
cycles $b_1,b_2$ as the transversal curves crossing the cycles
$a_1,a_2$ only. Our lemma is proved.

 Therefore the theorem is also proved.

{\bf Example: The case $g=3$}. A family of concrete foliations of the class
$T^2$ were demonstrated above for the real nonsingular algebraic
curves $$w^2=\prod (z-z_1)\ldots (z-z_8)=R(z), z_j\in R$$ with
the cycles $$a_1=[z_2z_3], a_2=[z_7z_8], a_3=[z_4z_5],
b_1=[z_1z_2], b_2=[z_6z_7]$$ and $\omega=P_2(z)dz/R(z)^{1/2}$ The
polynomial $P_2(z)=u+iv$ should be chosen such that $u\neq 0$
 in the segments $[z_2z_3], [z_4z_5],
[z_6z_7]$, and $v\neq 0$ in the
segments $[z_1z_2],[z_7z_8]$.

{\bf Question:} How to extend this basis to the
Transversal Canonical Basis?

 After cutting the Riemann Surface along the cycles $a_j,b_q$ we
realize it as a plane domain with following components of the
boundary:

 The external boundary $$A_1= (a_1b_1^{-1})\bigcup
(a_1^{-1}b_1)=A_1^+\bigcup A_1^-$$

 The internal boundary
$$A_2=(a_2b^{-1}_2)\bigcup (a_2^{-1}b_2)=A_2^+\bigcup A^-_2$$

  The
interim boundary $a_3^+\bigcup a_3^-$ inside.

 Our notations are  such that trajectories  enter this domain through the pieces with
sign $+$ and leave it through the pieces with the sign $-$.  We make a numeration
such that: $$2|A_1^{\pm}|=A_1> A_2=2|A_2^{\pm}|,
|a_3|=|a_3^{\pm}|=a$$ Here $A_k$ means also the measure of this
boundary component.  Nearby of the ends of the segments $A_l^+$
the trajectories enter our domain and almost immediately leave it
through the piece $A_l^-$. Therefore, there exist the first points
in $A^{\pm}_l$ where this picture ends. These are the endpoints
$x_{1,j}^{\pm},x_{2,j}^{\pm}\in A_j^{\pm}$ of the pair of
separatrices of  saddles (see Fig 5).

A complete list of  topologically  different types of the Plane
Diagrams in the class $T^2$ for the genus $g=3$ can be presented.
It shows that there is only one type such that  we cannot extend
the incomplete transversal basis to the complete transversal basis
(see Fig 6,a)): no closed transversal curve exists in this case
crossing the cycle $a_3$ and not crossing other curves
$a_1,b_1,a_2,b_2$ (i.e. joining $a^+$ and $a^-$ on the plane
diagram). In all other cases such transversal curve $b_3$ can be
constructed.

 Consider now this special case  where the
transversal incomplete basis cannot be extended (Fig 6,a).
\begin{lem}
There exists a  reconstruction of this basis such that the new
basis can be extended to the complete transversal canonical basis

\end{lem}

For the proof of this lemma, we construct
 a closed transversal curve $\gamma$  such that
it crosses the cycle $a=a_3$ in one point and  crosses  the
segment $A^+$  leaving our plane diagram (see Fig 7). It enters
$A^-$ in the equivalent point, say, through the cycle $b_1$. We
take a new incomplete transversal basis $a_2,b_2, a_3, \gamma,
a_1$. If $\gamma$ crosses $A^{\pm}_1$ through the cycle $a_1$, we
replace $a_1$ by the cycle $b_1$ as a last cycle in the new
incomplete basis of the type $T^2$. The proof follows from the
plane diagram of the new incomplete basis (let us drop here these
technical details).

 Comparing these lemmas with the construction
of special foliations in the previous section on the real Riemann
Surface, we are coming to the following

{\bf Conclusion.} For every real hyperelliptic Riemann Surface of
the form $w^2=\prod (z-z_1)\ldots (z-z_8)=R(z), z_k\neq z_l\in R$,
every generic form $\omega= P_2(z)dz/R(z)$ defines foliation
$\omega^R=0$ of the class $T$ if real and imaginary parts of the
polynomial $P=u+iv$ do not have zeroes on the cycles indicated
above for   genus $g=3$.

Extending these arguments, the present author proved
the existence of Transversal Canonical Basis for $g=3$ (see
 Appendix in \cite{NOV},
version 1). Very soon G.Levitt informed us that he can
 prove the existence of TCB
for all $g$ (see Appendix in \cite{NOV}, versions 2 and 3).
 We present here this proof based on the result of  \cite{L}.

Consider any  $C^1$-smooth vector field on the compact smooth
Riemann Surface $M_g$ of the genus $g\geq 1$, with saddle
(nondegenerate) critical points only. We introduce  {\bf the class
$G$ of vector fields} requiring that there is no saddle
connections and no periodic trajectories for the vector fields in
this class. Let a finite family of smooth disjoint

\begin{thm} For every generic dynamical system
of the class $G$ there exist  Transversal Canonical Basis
\end{thm}

 Proof. The result of \cite{L} is following: For every dynamical system of the class $G$
there exist exactly $3g-3$ non-selfintersecting and pairwise
non-intersecting  transversal cycles $A_i,B_i,C_i, i=1,...,g-1,$
where $C_{g-1}=C_1$. These cycle bound two sets of surfaces:
$A_i\bigcup B_i\bigcup C_i=\partial P_i$ where $P_i$ is a genus
zero surface (''pants''). The trajectories enter pants $P_i$
through the cycles $A_i,B_i$ and leave it through the cycle $C_i$.
There is exactly one saddle inside of $P_i$ for each number $i$.
The cycles $A_i,B_i,C_{i-1}$ bound also another set of pants $Q_i$
where $\partial Q_i=A_i,B_i,C_{i-1}$. The trajectories enter $Q_i$
through $C_{i-1}$ and leave it through $A_i$ and $B_i$. There is
also exactly one saddle inside $Q_i$.

The construction of Transversal Canonical Basis based on that
result is following: Define  $a_g$ as $a_g=C_1$. Choose a segment
of trajectory $\gamma$ starting and ending in $C_1$ assuming that
it passed all this ''necklace'' through $P_k,Q_k$. For each $i$
this segment meets either $A_i$ or $B_i$. Let it meets $B_i$ (it
does not matter). We define $a_i$ as $a_i=A_i,i=1,...,g-1$. Now we
define $b_g$ as a piece of trajectory $\gamma$ properly closed
around the cycle $C_1=a_g$. This curve is almost transversal in
our terminology (above), so its natural small approximation is
transversal. We are going to construct the cycles $b_i$ for
$i=1,...,g-1$ in the union $P_i\bigcup Q_i$. Consider the saddle
$q_i\in Q_i$. There are two separatrices leaving $q_i$ and coming
to $A_i$ and $B_i$ correspondingly. Continue them until they reach
$C_i$ through $P_i$. Denote these pieces of separatrices by
$\gamma_{1,i},\gamma_{2,i}$. Find the segment $S_i$ on the cycle
$C_i$ not crossing the curve $\gamma$ chosen above for the
construction of the cycle $b_g$. The curve
$\gamma_{1,i}S_i\gamma_{2,i}^{-1}$ is closed and transversal
everywhere except the saddle point $q_i$. We approximate now the
separatrices $\gamma_{1,i}$ and $\gamma_{2,i}^{-1}$ by the two
pieces of nonseparatrix trajectories $\gamma'_i,\gamma''_i$
starting nearby of the saddle $q_i$ from one side of the curve
$\gamma_{1,i}\gamma_{2,i}^{-1}$. Close this pair by the small
transversal piece $s_i$. There are two possibilities here (two
sides). We choose the side such that the orientation of the
transversal piece $s_i$ near the saddle $q_i$ is agreed with the
orientation of the segment $S_i$, so the whole curve $b_i=s_i
\gamma'_iS_i\gamma''_i$ is closed and almost transversal. We
define $b_i$  as a proper small transversal approximation of that
curve. So our theorem is proved because every cycle $b_i$ crosses
exactly one cycle $a_i$ for $i=1,2,...,g$.

\vspace{0.3cm}

\centerline{\Large Section 3. Topological Study of the Class $T$.}

\vspace{0.2cm}

Let us describe here some simple general topological properties of
the class $T$ foliations
 After cutting the surface $V$ along the transversal
canonical basis $a_j,b_j$, we are coming to the fundamental domain
$\bar{V}$ of the group $Z^{2g}$ acting on the maximal abelian
covering $V'$. There is exactly $g$ boundary ''squares'' $\partial
\bar{V}=\bigcup_{j} A_j$ where $A_j= a_jb^{-1}_ja_j^{-1}b_j
=A^+_j\bigcup A^-_j$. Every piece $A_j^{\pm}$ consists of  exactly
two basic cycles $a_j,b_j$ attached to each other. These pieces
are chosen such that trajectories enter $A^+_j$ from outside
through $a_jb_j^{-1}$ and enter the fundamental domain
$\bar{V}$.  In the areas nearby of the ends the trajectories
enter $A^-_j$ from inside and leave fundamental domain. There is
exactly $2g-2$ saddle points inside of $\bar{V}$. They are no saddles
located on the selected transversal cycles $a_j,b_j$. Our
foliation nearby of each boundary square $A^+$ looks exactly as a
straight line flow. It means in particular that there exist two
pairs of points $x_{1,j}^{\pm},x^{\pm}_{2,j}\subset A^{\pm}_j$
which are the endpoints of separatrices in $\bar{V}$, nearest to
the ends at the each side $A^{\pm}_j$ (see Fig 3 and 5). We call
them {\bf The Boundary Separatrices} belonging to {\bf The
Boundary Saddles } $S_{j,1},S_{j,2}$ for the cycle $A_j$. We can see
$2g$  of such saddles $S_{j,1},S_{j,2}$ looking from the
endpoints  of $A_j^{\pm}$,  but some of them are in fact the same. At
least two of them should coincide. We call corresponding saddles {\bf
The ''Double-boundary'' saddles}.

{\bf Definitions}.

 1.The saddle point $S \in \bar{V}$
has a type $<jklm>$ if it has two incoming separatrices starting
in   $A^+_j,A^+_l$ and two outcoming separatrices ending in
$A^-_k,A^-_m$. The indices are written here in the cyclic order
corresponding to the orientation of $\bar{V}$. Any cyclic
permutation of indices defines an equivalent type. we normally
write  indices   starting from the incoming separatrix,
 $<jklm>$ or $<lmjk>$.

 2. We call
foliation {\bf Minimal} if all saddle points in $\bar{V}$ are
  of the boundary saddles. In particular, their
types are $<jjkl>$. We call foliation {\bf  Simple} if there are
exactly two saddle points of the double-boundary  types $jjkk$ and
 $jjll$ correspondingly. The index $j$ we call {\bf Selected}. We
say that foliation has a rank equal to $r$, if there exists
exactly $r$ saddles of the types $<jklm>$ where all four indices
are non-boundary. In particular, we have $0\leq r\leq g-2$. There
is exactly $ t$ saddles of the double-boundary types
where $t-r=2$. There is also $2g-2t$ other saddles of the types
like $<jjkl>$ where only index $j$ corresponds to the boundary
separatrices. The extreme cases are $r=0,t=2$ which we call {\bf
special} (above), and $r=g-2, t=g$ which we call {\bf maximal}. In
the maximal case there exists a maximal number of saddles of
 the nonboundary types, and
all boundary type saddles are organized in the pairs. One might
say that for the maximal type every index is selected. {\bf For
maximal type the genus should be an even number} because the
boundaries $A^j$ are organized in the cycles. Every cycle
should contain even number of them, by the elementary orientation
argument. The relation $2g-t+r=2g-2$ for the total number of
saddles gives $t-r=2$. For the case $g=2$ we obviously have $r=0$.
For the case $g=3$ the only possible case is $t=2,r=0$ (the
simple foliations); the case $t=3,r=1$  cannot be realized for
$g=3$ because it is maximal in this case but the  maximal
case corresponds to the even genus only.

\vspace{0.2cm}

{\bf How to Build Hamiltonian Systems from the elementary pieces.}\\
Let us introduce  following {\bf
Building Data} (see Fig 9):

I.{\bf The Plane Diagram} consisting of  the generic Hamiltonian
System on the 2-sphere $S^2$ generated by the hamiltonian $h$ with
nondegenerate critical points only (centers and saddles) sitting
on different levels (see fig 9a,9b,9c). Let minimum of $h$  is located in the point $0$,
and maximum in $\infty$. There are $t-2$ other centers and $r=t-2$ saddles.
Let exactly $g$  oriented segments are given transversal to the system
everywhere except centers
$t_1,...,t_g \subset S^2$ with  measures
$m_1,m_2,...,m_g$ provided by hamiltonian $h$, such that:

a.They do not cross each other;  the values of Hamiltonian in all
their ends, all centers and all saddles  are distinct except that exactly
 two of them
meet each other in every center; They do not touch any saddle
 on the two-sphere.

b.Every cyclic and separatrix trajectory of the hamiltonian system
on $S^2$
 meets at least one of these segments.

 We make cuts along these segments and define the sides $t_j^{\pm}$
 where trajectories  leave and enter them correspondingly.

II.{\bf The Torical  Data} consisting of the $g$ tori
$T^2_j,j=1,...,g$, with  irrational straight
line flows
 and selected oriented transversal
segments $s_1,...,s_g$ (one for each torus). Their  transversal
measures are equal to the same numbers $m_1,...,m_g$. The
Transversal Canonical Basis  $a_k,b_k$ in every $T^*_{m_k}$ equal to torus minus
the segment $s_k$, is selected. Here
 $a_k$ are positive, $b_k$ are negative, and
$|a_k|+|b_k|>m_k$.
We make  define
the sides $s^{\pm}_k$ in the boundary $\partial T^*_k$  as for the segments $t_j$.

Identifying the segments $s^+_j$ on the tori with $t^-_j$ on the
sphere $S^2$ and vice versa, we obtain a Riemann surface $V=M^2_g$ with
Hamiltonian System (Foliation).
We glue them such that every end of the
segment $t_j\subset S^2$ located on the trajectory $\gamma\subset S^2$,
defines exactly one saddle of the
boundary type $<jjkl>$. Here $t_k,t_l$ are  the  segments where this piece of
the  trajectory $\gamma$ ends  on the
2-sphere $S^2$. Two pieces of $\gamma$ divided by the segment $t_j$, provide  a pair of
non-boundary separatrices for the saddle on the Riemann Surface defined by this picture.

By construction, every center  generates
 a double-boundary
saddle of the type $<jjkk>$. So we have $t-r=2$.

 Let us perform following operations:

1.Cut our surface along the TCB. The boundary of the domain
$\bar{V}$  is equal to the union $\bigcup_j (A^+_j\bigcup
A^-_j)=\partial \bar{V}$. Every component is presented as boundary of
the Fundamental Parallelogram $P_j$ associated with the 2-torus $T^2_j$. Our
flow covers the boundary  of $P_j$ as a straight line flow: the
trajectories enter through the  pair of cuts
$a_jb_j^{-1}$ and leave through the pair $a_j^{-1}b^j$ (see Fig
2).

2.Find for every 2-torus $T^2_j$ the pair of  boundary saddles in
$P_j$ and join them by the pair of transversal segments $s_j^{\pm}\subset P_j
$ along the parts $A^{\pm}_j$ (see
Fig 5). We perform this operation in the fundamental parallelogram
$P_j$. These segments should be chosen in
such a way that outside of them in $P_j$  near the one-skeleton, we have a
straight line flows.

3.Cut our surface $V=M^2_g$ along  the segments $s^+_j\bigcup
s^-_j$.  It is divided now into the torical pieces $T^2_j$ and one
 piece $S^2$ whose boundary consists of the curves
  $\bigcup_j(t^+_j\bigcup t^-_j)$  for $S^2$, and $s^+_j\bigcup s^-_j$
 for the tori $T^2_j$. After cutting the surface $M^2_g$ along the
 pieces $s_j^{\pm}$ we keep the notation $s^{\pm}_j$ for the curves
 in the tori, but for the  part $S^2$ we
 denote them by $t^{\pm}_j$.

4.Glue $t^+_j$ with $t^-_j$ for the sphere $S^2_j$,  and
$s^+_j$ with $s^-_j$ for the tori $T^2_j$, preserving the
transversal measure.  The system on
 2-sphere appears with $g$ selected transversal segments $t_j$ and
  2-tori $T^2_j$ with  straight line flows where the
transversal pieces $s_j$ with $m_j$ are marked.

5.Near the double-boundary saddle  this picture is
topologically equivalent to the center, but this equivalence in
non-smooth.

 We can see  that our construction  preserves  the measure-type
invariants.

Therefore we are coming to the following

\begin{thm} Every generic foliation given by the real part of
 holomorphic one-form,
can be obtained by the measure-preserving gluing of the
 pieces $(S^2,h,t_1,...t_g)$
and $(T^2_1,s_1),...,(T^2_g,s_g)$ along the transversal segments
$t_j$ and $s_j$, as it was described above. For the genus $g=2$ we
can remove  sphere $S^2$ from the description: Every generic
foliation given by the real part of holomorphic one-form, can be
obtained from the pair of tori $(T^2_1,s_1)$ and $(T^2_2,s_2)$
with different irrational straight line flows and  transversal
 segments $s_1,s_2$ with transversal measure $m_1=m_2=m$
\end{thm}

{\bf Example. The Topological Types of the Minimal Foliations}.

 For
the Minimal Foliations above we have $t=2,r=0$. The hamiltonian
system on the 2-sphere is trivial (see Fig 9,a): It can be
realized by the rotations around the points $0$ and $\infty$. There are no
saddles on the sphere. For the simplest case $g=2$ we have two segments
$t_1,t_2$. Both of them join $0$ and $\infty$. So they form a
cycle of the length 2. The difference between the values of
Hamiltonian $h(\infty)-h(0)$ is equal to $m_1=m_2$. All possible
pictures of the transversal segments can be easily classified here
for every genus $g$ (see Fig 9,a and b).

There exist following types of topological configurations only:

a.The Plane Diagram has exactly one ''cycle'' $t_1,t_2$ of the
length two (like for $g=2$) and $g-2$ disjoint segments $t_j,j\geq
3$; This type is available for all $g\geq 2$ (see fig 9a for $g=2$).

b.The Plane Diagram has two pairs $t_1,t_2$ and $t_3,t_4$ where
the members of each pair meet each other either  in the center $0$
or in $\infty$. Other $g-4$  segments $t_j,j\geq 5$ are disjoint. In the
second case we have $g\geq 4$.

c.The Plane Diagram has exactly one connected set consisting of 3
segments $t_1t_2t_3$ passing through $0,\infty$. Other $g-3$
 segments $t_j,j\geq 4$ are disjoint. For this type we have $g\geq 3$.

\begin{thm} For $g=2,3$ every class $T$ foliation is simple.
The Maximal
type exists only for even genus $g\geq 4$
\end{thm}

Proof. For $g=2$ it was already established
above: all generic (irreducible) hamiltonian foliations are
simple. Both indices $k=1,2$ are selected. Consider now the case $g=3$. Our foliation can
be either simple ($t=2,r=0$) or maximal ($t=3,r=1$) in this case.
We have $t=3,2g-2t=0$ for the  maximal case. If it is so, every
boundary saddle should be paired with some other. So there is a
cyclic sequence of boundary saddles containing all three boundary
components. However, every cyclic sequence
 should contain even number of boundary components ( and the same
number of boundary saddles), otherwise the orientation of
foliation is destroyed. This is possible only for even number of
indices which is equal to genus. Our conclusion is that $g\geq 4$.
This theorem is proved.

The  2-sphere $S^2$ is covered by the non-extendable {\bf '' corridors''}
between two  segments $t_j,t_k\subset S^2$ (see fig 9b) consisting
of the non-separatrix trajectories  moving from the inner
points of $t_j$ to the inner points of $t_k$ not meeting  any
points of $t_l$  and the saddles inside. The right and left sides of
 corridors are either separatrices of the saddles in $S^2$ or
the trajectories passing  through the ends of some
segments $t_l$.

Classification of the generic Morse functions $h$ on the sphere
$S^2$ is following: Take any connected
trivalent finite tree $R$ (see fig 9c). It has vertices divided into the  $r$
{\bf Inner Vertices} and $t$ {\bf Ends}. Assign to each vertex
$Q\in R$ a value $h(Q)$; these values   are not equal to each
other $h(Q)\neq h(P),P\neq Q$; therefore the edges become
oriented  looking ''up'' where  $h$ is increasing.
Every inner vertex $Q$ is a ''saddle'', i.e.
 $$\min_ih(Q_i)<h(Q)<
\max_ih(Q_i), i=1,2,3$$ where $Q_i$ are the neighbors. The
function $h$ on the graph should be monotonic  for every edge
$[Q_1Q_2]\subset \tilde{R}$.

For the description of the set of transversal segments $t_j\subset
S^2$ we introduce  {\bf locally constant set-valued function}
$\Psi (P),P\in \tilde{R} $ on the graph $\tilde{R}$: its values $\Psi(P)$ are the
cyclically ordered nonempty finite sets with properties indicated below:
$$\Psi(P)=\{...t_1(P)<t_2(P)<...<t_q(P)<t_1(P)...\}$$
The number of points $q=|\Psi(P)|$ is equal
to two nearby of the centers $Q_j$ collapsing to one in the endpoint $Q_j$
 The number $q$ may change
$q\rightarrow q\pm 1$ in the isolated points $P_l\in \tilde{R}$ such that
$h(P_l)\neq h(Q)$ for all inner vertices $Q$. Passing ''up''
through the inner vertex $Q$ this set either splits into the pair
of cyclically ordered sets $\Psi(P)\rightarrow \Psi_1\bigcup
\Psi_2$ inheriting orders  where $|\Psi_1|+|\Psi_2|=|\Psi(P)|$, or
some pair of cyclically ordered sets is unified into the one
ordered set choosing some initial points in each of them (the
inverse process). Every continuous ''one-point branch'' $t_j(P)$
living  between the points $P_1,P_2\in \tilde{R}$,
defines the transversal segment $t_j$. Its transversal measure is
equal to $m_j=h(P_2)-h(P_1)>0$.\\
 {\bf Example: The Maximal Foliations} for $g=4$.\\
For the maximal  foliations of the class $T$ all boundary saddles
are paired with each other. It means that all of them organize a
system of cycles where the next boundary saddle is paired with the
previous one. The length of every cycle  is an even number
 $2l_q,q=1,\ldots,f$, so $2l_1+\ldots +2l_f=g$. We say that
  system has a cycle type $(l_1,\ldots,l_f)$. The maximal
system also
 $g-2$ saddles  such that all 4 separatrices are nonboundary.
 For the case $g=4$ we have two possibilities for the
cycles, namely $(1,1)$ and $(2)$. The type $(2)$ is especially
interesting.  This cycle separates a 2-sphere on the South and
North Hemisphere (see Fig 12) where $A_1,A_2,A_3,A_4$ are located
along the equator. The additional two saddles are sitting in the
poles  with separatrix curves going to the each ''country''
$A_k,k=1,2,3,4$ along the 4 selected meridians from the north and
south poles.

\begin{lem} The  transitions $A^+_k\rightarrow A^-_l$
with measures $m_{kl}$
 visible from the poles (i.e. located in the corresponding hemisphere)
 are the following:

From the North Pole we can see the transitions $$A_k^+\rightarrow
A^-_l,k=1,3,l=2,4$$
 From the South Pole we can see the transitions
$$A^+_k\rightarrow A^-_l,k=2,4,l=1,3$$ They satisfy to the
Conservation Law $\sum_km_{kl}=|A_l|,\sum_lm_{kl}=|A_k|$, where
$k$ and $l$ are neighbors in the cyclic order $...12341...$.
We have $A_1-A_2+A_3-A_4=0$. The constant
flux in one direction of the cycle $1234$ is defined provided by
the asymmetry $m$ of the  measures $m_{kl}-m_{lk}$:
$$m_{12}-m_{21}=m_{23}-m_{32}=m_{34}-m_{43}=m_{41}-m_{14}=m$$
We say that the  system is rotating clockwise if $m>0$. It is
rotating contrclockwise  if $m<0$.
\end{lem}

\vspace{0.2cm}

\centerline{\Large Section 4.Three-Street Picture on the
torus with  obstacle.}

\vspace{0.2cm}

We are dealing now with one torus $T^2=T^2_k$ with generic straight line flow.
 The torus is presented as a lattice in $C=C_k$ given by the parallelogram $P=P_k$,
 the flow is a set of vertical lines. A transversal segment $s=s_k$ with
  measure $m=m_k$ is marked in $T^2$.\\
{\bf Question}. How long the trajectory can move in the plane $C$ (i.e. in
the torus $T^2$) until it hits some periodically repeated copy
of the  segment $s$ ?

Every trajectory starts and ends in some segments of the selected periodic
family generated by the segment $s$. Family of such pieces
with fixed ends  form the connected strip. We require that these
strips cannot be extended to the left and right: every trajectory
in the strip ends in the same  segments. The right and the left boundaries of the strip
are separatrices of some saddles.  Every such
strip has  {\bf Width $w$}. It is
equal to the transversal measure.

 {\bf Definition}.We call the unextendable strips by {\bf The Streets}
 and denote them $p^{\tau},\tau=0,1,2$. We  denote a longest unextendable strip
by  $p^0$. Its right and left boundaries   meet the ends of some segments $s',s''$
 of our family
strictly inside (see Fig 11). The
 lower and upper segments of this strip should be located
strictly inside of the corresponding segments $s, s'''$. Their
 Widths we denote by $ |p^{\tau}|$
correspondingly. The street number 2 is located from the right
side from the longest one, the street number 1--from the left
side. We assign also homology classes to the street $h^{\tau}\in H_1(T^2,Z)$
naturally.

\begin{lem}
For every straight-line foliation of $T^*_m$, i.e. of the torus $T^2$ with transversal
canonical basis $a,b$ and obstacle $s\subset T^2$ with measure $m<|a|+|b|$,
there exists exactly three streets
$p^{\tau},\tau=0,1,2$ such that $\sum_{\tau } |p^{\tau}|=m$
and $h^1+h^2=h^0$. Two smaller streets are attached to the
longest one from the right and left sides. This picture is
invariant under the involution changing time and orientation of
the transversal segments. The union of the three streets started
in the segment $s$ is a fundamental domain of the
 group $Z^2$ in the plane $C$ (see Fig 11).
\end{lem}
Proof. Construct first the longest street $p^0$. We start from
any non-separatrix trajectory ending in  segments of our family
inside both
of them. Extend this strip to the left and right) directions.
We either meet the ends of upper or lower segments or the end of some
segment  strictly between. In the first case we see
that extending beyond the  left end we  construct longer strips.
  start the same process with longer strip.
Finally, we reach the locally maximal strip $ss'''$.
We
 meet from the both sides of it the ends of some segments $s',s''$,
otherwise its height would  be extendable. So the
maximal street is constructed. Consider the neighboring streets
from the right and left sides.
They are shorter. Denote the left one by $p^1$ and the
right one by $p^2$. They
 can be extended to the  left and right boundaries of the
 locally longest one. This statement follows from the
periodicity of the system of segments.
 The same argument can be applied to the upper end $s'''$.
 We  see that our locally longest
street is surrounded by the exactly  four unextendable streets
(two streets from each side). They are restricted from the right
and left sides by the ends of the lower or upper segments $s,s'''$(see Fig
11 ). The pairs of streets located across the diagonal of each
other are equal.
So, there are no other unextendable streets except these three.

Lemma is proved.

\vspace{0.2cm}

\centerline{\Large Section 4. The Homology and Homotopy Classes.}

\vspace{0.3cm}

How to describe  the image of Fundamenal Transversal Semigroups in
fundamental group $\pi_1^+(\Re)\rightarrow \pi_1(V)$ and in
homology group $H_1(V,Z)$?

There are three types of the transversal curves generating all of
them:
(I) The ''Torical Type'' transversal curves $\gamma\subset T^*_m\subset V$
in one 2-torus not touching   the
obstacle $s$ of the transversal measure $m$;
(II) The ''Trajectory Type'' transversal curves $\gamma\in V$ (or the Poincare
Curves). They coincide with some trajectory of the Hamiltonian
System started and ended in the transversal interval $s^+_k$ in the
 the torus $T^2_k$. We make it closed
joining the endpoints by the
shortest transversal interval along $s^+_k$;
(III) The general   non-selfintersecting transversal closed
curves.

{\bf Step 1}: Classify all closed transversal curves $\gamma\subset T^*_m$.

Let $|a|,|b|$ denote the measures of the basic transversal
closed curves $a,b$. {\bf We made our  notations
above such that the cycles $a, b^{-1}$ are  positive and
transversal.} According to our construction, the measure
$m$ of the segments $s$ satisfies to the inequality
$0<m<|a|+|b|$. Define  minimal  nonnegative
pairs of integers $(u>0,v\geq 0), (w\geq 0,y> 0)$ such
that
$$m>u|a|-v|b|>0,\\\\m>y|b|-w|a|>0$$ (see Fig
13). It means exactly that  these new lattice vectors represent
the shifts $s',s''$ of the segment $s$
 visible directly from $s$ along the shortest trajectories
 looking to the positive time
direction.
\begin{lem}
The  new lattice vectors $h^1,h^2$ have transversal measures
equal to $|a^*|=u|a|-v|b|$ and $|b^*|=y|b|-w|a|$ where $|a^*|+b^*|>m$.
 Their homology classes
$h^1=u[a]+v[b]$ and $-h^2=w[a]+y[b]$, generate homological image of
semigroup of positive closed transversal curves not crossing the
segments $s$ (i.e. of the Torical Type (I)). These classes  have
canonical lifts $a^*,b^*$ to the free group $F_2$ generated by
$a,b$; The lifts $a^*,(b^*)^{-1}$  generate free semigroup of all homotopy
classes of positive closed transversal curves not crossing the
segments $s_k^{\pm}$, starting and ending in the street number
$0$. These semigroups depend on the triple $|a|,|b|,m$
 only.
 Similar classes $h^1_-,h^2_-$ can be constructed for the negative
 time direction: their homology classes  are opposite
to the positive ones.
 They represent
the same streets going back. Their lifts $a'_-,b'_-$ to
the free fundamental group $\pi_1(T^*_m,s^-)$ of torus with cut
 along $s$,
are defined as a mirror
symmetry of the lifts $a',b'$. Here
$t \rightarrow -t, s^+\rightarrow -s^-$.
\end{lem}

Proof. We have the streets $p^{\alpha}$ in the
plane $C$  where the longest
one is $p^0$. It is located between two others. The right one is
$p^2$, and the left one is $p^1$. For the homology classes we have $h^1+h^2=h^0$,
 and for the  transversal
measures  $\sum_{\alpha}|p^{\alpha}|=m$. Their bottom
parts cover together the segment $s$. This picture is invariant
under the change of direction of time and simultaneous permutation
of the lower and upper segments $s$. All positive transversal paths  in the
plane $C_k$ can be written combinatorially, numerating the streets which they cross,
 in the form $$\ldots
\rightarrow\alpha_q \rightarrow 0\rightarrow
\alpha_{q+1}\rightarrow 0\rightarrow\alpha_{q+2}\rightarrow
0\rightarrow \ldots$$ where $\alpha_q=1,2$. We have a pair of
positive ''basic cycles'' $[0\rightarrow 1\rightarrow
0]=(b^*)^{-1}$ and $[0\rightarrow 2\rightarrow
0]=a^*$. All other transversal cycles not touching the
segment $s$ and its shifts, starting and ending in the longest
strip $0$, have a form of the arbitrary word in the free semigroup
generated by $a^*,(b^*)^{-1}$ in the
plane $C$. The measures of these new $m$-dependent basic cycles
are
$$|a^*|=|p^0|+|p^2|,|b^*|=|p^0|+|p^1|$$ Topologically the cycles
$a^*,b^*$ represent some canonical
 $m$-dependent basis of
 elementary shifts in the group $Z^2$. To define the elements $a^*,b^*$,
 lifting them to fundamental group,
 we choose the simplest transversal
paths joining the initial point in $s^+$ with its image in $C$ avoiding obstacles.
 New basic shifts map the segment
$s\subset P_{0,0}$ exactly into the the segments $s',s''$
attached to the middle part of the long street from the left side
(for $a^*$) and from the right side (for $b^*$)--see Fig 13.
 Homologically  the
elements $[a^*]=h^1,[b^*]=h^2_k$
 are calculated in
the formulation of lemma using the transversal measures
$|a|,|b|$ of the initial basic cycles and the measure $m$. From geometric
description of new paths $a^*,b^*$  we can see
that they satisfy to the same relation as the original $a,b$:
Their commutator path exactly surrounds one segment $s$ on the
plane. Therefore our lemma is proved.

This is the end of the Step 1.

Our 2-torus with cut along the obstacle is a part of the Riemann surface
 $T^*_{m_k}\subset V$. As we know, in order to have  well defined
element $\{\gamma\}\in\pi_1(V)$, we need to have both ends belonging to $s^+_k$.
 To
define homology class $[\gamma]$ is enough to have both ends belonging
to $s^+_k\bigcup s^-_k$ because
 $H_1(V,s^+_k\bigcup s^-_k)=H_1(V)$. Therefore we arrive to the following\\
{\bf Conclusion}: {\it Every trajectory of  Hamiltonian System
on the Riemann Surface $V$
constructed out of  Building Data above, defines an infinite ''sum'' of homology classes
$(...+[\gamma_N]+[\gamma_{N+1}]+[\gamma_{N+2}]+...)$ where
$[\gamma_N]\in H_1(V)=H_1(V,S^2)=\sum_k H_1(T^2_k)$. Every class $[\gamma_N]$ is equal
to the homology class of the street $h^{\tau}_k\in H_1(T^2_k)$ for $\tau=0,1,2$. So we have
 the exact coding of trajectories by the sequences of symbols $\tau,k$ if these homology
 classes are known.}

\vspace{0.2cm}

{\Large{\bf The case $g=2$}}. For the special case $g=2$  we can assign the invariant
$\phi_{\alpha\beta}\in\pi_1(V,s_1^+)$ to every piece of trajectory
passing two streets $$\gamma_q\subset
p^{\alpha}_1p^{\beta}_2,\phi:\gamma_q\rightarrow
\pi_1(V,s^+_1)$$ So every infinite trajectory $\gamma$ can be viewed as an
 infinite product of pieces $$\gamma=\ldots
\phi_{\alpha_q\beta_q}\phi_{\alpha_{q+1}\beta_{q+1}}\ldots$$
 Every finite connected  piece $\gamma'$ of this
sequence defines the element $\phi(\gamma')\in \pi_1(V,s^+_1)$. Its
ends can be joined by the shortest transversal piece along the
segment $s$. Depending on orientation of this piece, either
$\phi(\gamma')$ or $\phi(\gamma')^{-1}\in\pi_1(V,s^+_1)$ defines a
closed positive transversal curve of the Poincare Type in the Fundamental Transversal Semigroup.

{\bf Reduction to the Standard Model}:  The segment $s$ of the
length $m$ is divided on  $5$ connected open pieces $\tau_q$:
There exist
 exactly 9 possible types of trajectory pieces
$\{\alpha\beta\},\alpha,\beta=1,0,2$ with measures
$p_{\alpha\beta}$ but 4 of them are in fact empty. They   have
 measure equal to zero. In order to see that, we remind how these pieces
 were constructed. A segment $s=s_1$ of the total measure $m$ is
 divided into 3 pieces by the points $0,1,2,3$ for $k=1$ and by
 the points $0',1',2',3'$ for $k=2, s=s_2$: For the streets  we have:
 $$ \alpha=1,0,2=[01],[12],[23],\beta=1',0',2'=[0'1'],[1'2'],[2'3']$$
 So,  the index $\alpha=1,0,2$ corresponds
 to the segments $[01],[12],[23] $ with measures $p^{\alpha}_1$,
 and the index $\beta=1',0',2'$  corresponds
 to the segments $[0'1'],[1'2'],[2'3']$ with measures $p^{\beta}_2$.
The positions of the points $0=0',3=3'$ are fixed. Other points
never coincide for the generic foliations.

 Every jump from $C_1$ to $C_2$  is accompanied by the permutation
of 3 segments: the left street number 1 in $C_1$ ends up in the
extreme right part of $s=s'$ before making jump to $C_2$. The
right street number 2 in $C_1$ ends up in the extreme left part of
$s=s''$ in $C_2$. So jumping from $C_1$ to $C_2$, we should
permute the segments $1=[01]=1$ and $2=[23]=2$ preserving
orientation:
$$\eta_{12}:2\rightarrow 0=2^*,
 3\rightarrow 3^*\in s,
 1\rightarrow 3=1^*, 0\rightarrow 0^*\in s$$.
 Here $ |[23]|=|[2^*3^*]|,|[01]|=|[0^*1^*]|$.
The segment $s$ is divided by the points
 $0=2^*,1',2',3^*,0^*,3=1^*$, so it is presented as a    union of   5 sub-segments
$$s=\tau_1+...+\tau_5$$
 In order to return
back from the second torus (plane) $C_2$ back to $C_1$, we need to
apply the similar map $\eta_{21}$ based on the permutation of the
streets $2'=[2'3']$ and $1'=[0'1']$. Our broken isometry
$i_{\sigma}$ based on the permutation $\sigma$  of 5 pieces, is
defined as a composition
$$i_{\sigma}=\eta_{21}\eta_{12}$$
 \begin{lem} There exist six Topological Types of Foliations
  (the measures $p_{\alpha\beta}$ of the
sub-segments $\tau_q,q=1,2,3,4,5$ are given
 in the natural order on the segment $s$):
$$(I):0=2^*<1'<2'<3^*<0^*<1^*=3 ; \sigma=(32541);$$
$$p_{12}+p_{02}+p_{21}+p_{20}+p_{22}=m$$
$$(II):0=2^*<1'<3^*<2'<0^*<1^*=3; \sigma=(24153);$$
$$p_{12}+p_{01}+p_{02}+p_{21}+p_{20}=m$$
$$(III):0=2^*<1'<3^*<0^*<2'<1^*=3;\sigma=(41523);$$
$$p_{10}+p_{12}+p_{00}+p_{21}+p_{20}=m$$
$$(IV):0=2^*<3^*<1'<2'<0^*<1^*=3;\sigma=(25314);$$
$$p_{12}+p_{01}+p_{00}+p_{02}+p_{21}=m$$
$$(V):0=2^*<3^*<1'<0^*<2'<1^*=q3;\sigma=(31524);$$
$$p_{10}+p_{21}+p_{01}+p_{00}+p_{21}=m$$
$$(VI):0=2^*<3^*<0^*<1'<2'<1^*=3;\sigma=(52134);$$
$$p_{11}+p_{10}+p_{12}+p_{01}+p_{02}=m$$ All other measures
$p_{\alpha\beta}$ are equal to zero. We have
$$\sum_{\alpha}p_{\alpha\beta}=p_2^{\beta};
\sum_{\beta}p_{\alpha\beta}=p_1^{\alpha}$$
\end{lem}

 As we can see, our 3-street model automatically
creates the standard type broken isometry
$\eta_{21}\eta_{12}=i_{\sigma}:s\rightarrow s$ generated by the
permutation $\sigma$ of 5 pieces $\tau_q$ of the segment $s$.
 Such systems were studied by the ergodic people since 1970s.  Let
us point out that our combinatorial model easily provides  full
information about the topology lying behind this permutation. We
know geometry of all pieces.

Let us define an algebraic object (no doubt,  considered by the
ergodic people many years ago): {\bf The Associative Semigroup}
$S_{\sigma,\tau}$ with ''measure''. It is generated by the 5
generators $R_1,...,R_5\in S_{\sigma,\tau}$. There is also a zero
element $0\in S_{\sigma,\tau}$. We define multiplication in the
semigroup $S_{\sigma,\tau}$ as in the free one but factorize by the ''zero
measure''  words: they are equal to zero. In order to
define which words are equal to zero, we assign to every generator
an interval $R_q\rightarrow \tau_q,q=1,...,5$.
 We assign to the word $R$ the set
$\tau_R$ where $\tau_{R_q}=\tau_q$ for generators. Let the word $R$ has a length $N$.
 By induction, we assign to
the word $R_pR$ the set $\tau_{R_pR}$ such that
$$\tau_{R_pR}=i_{\sigma}^{-N}(\tau_p)\bigcap \tau_R$$ if it is
nonempty. We put $R_pR=0$ otherwise.
This semigroup is associative. For every word in the
free semigroup $R=R_{q_1}...R_{q_N}$ we assigne a set
$$i_{\sigma}^{-N+1} (\tau_{q_1})\bigcap
i_{\sigma}^{-N+2}(\tau_{q_2})\bigcap...\bigcap \tau_{q_N}=\tau_R$$
whose measure is well defined. We put $R$ equal to zero if its
measure is equal to zero.

The ordered infinite sequence $R_{\infty}=\prod_{p\in Z} R_{q_p}$,
 defines trajectory  if and
only if  every finite sub-word is nonzero. The measure is defined
for the ''cylindrical'' sets $U_R$ consisting of all
''trajectories'' with the same sub-word $R$ sitting in the same
place for all of them. It is equal to
$\tau_R$. The ''Shift Function'' is well-defined by the
semigroup.
{\bf Our goal is to calculate representation of this
semigroup in the fundamental group of the Riemann Surface
generated by the Hamiltonian System.}

 Our global reduction of the flow
is based on the pair of non-closed transversal
 segments $s^{\pm}$ leading from one saddle to another.
This construction seems to be  the best possible genus 2 analog
 of the reduction of straight line flow on the 2-torus  to the
rotation of circle.  Many features of this construction certainly
appeared
 in some very specific examples
 studied before.

We always choose initial point in the segment $s^+_1\subset C-1$.
 It is equivalent   to the choice of initial point in the segment
$s^-_2=s^+_1\subset C_2$.  {\bf The New $m$-Dependent
Transversal Canonical Basis} is (see Fig 13):
 $$a^*_1,b^*_1,a^*_2=a^*_{2,-},b^*_2=b^*_{2,-}$$
  This basis is attached to
 the segment $s^+_1=s^-_2$, so we
treat them as the elements of fundamental group $\pi_1(V,s^+_1)$.
It was proved before that  homology classes of the  (positive) streets
  are following:
$$h^1_k,h^2_k,h^0_k=h^1_k+h^2_k\in
H_1(V,Z)$$ Here $k=1,2$. For the
negative time we simply change sign
$[a^*_-]=-h^k_1,[b^*_-]=-h^k_2$.

 Let us describe the homotopy classes
$\phi_{\alpha\beta}\in \pi_1(V,s^+_1)$
 of the two-street paths $\alpha\beta=p^{\alpha}_1p^{\beta}_2$ with ends in the same open segment
$s=s^+_1=s^-_2 \in T^2_1\bigcap T^2_2$.

 Consider first the plane $C_1\rightarrow T^2_1$.
 In order to make closed paths out of the streets $p^{\alpha}_1$,
 we join ends going
 around the segment $s$ from $s^-_1$ to $s^+_1$ from the
 right or from the left side of it (see Fig 14):

For the   street number 1 we construct  the  path
$a^*_1$ extending the street around $s'$ from
 the right side along the path $\kappa_1$, circling contr-clockwise around $s'$;
For the street number 2 we construct  the path $b^*_1$ closing street
around $s''$ from
 the left side along the path $\kappa_2$, circling clockwise around $s''$;
 So we have $p^1_1\sim a^*_1\kappa_1^{-1}$,
$p^2_1\sim b^*_1\kappa_2^{-1}$. Here and below {\bf  symbol
$\sim$ means ''homotopic with fixed ends''}.
 We define a closed path circling clockwise around the segment $s$ in $C_1$:
 $$\kappa=\kappa_2\kappa_1^{-1}\sim (s^+\bigcup s^-) \sim
a^*_1b^*_1(a^*_1)^{-1}(b^*_1)^{-1}$$ We
assign   element $p^0_1\sim a^*_1b^*_1\kappa_2^{-1}$ to the street number zero.
Same
description we have also for the streets
$p^{\alpha}_2(-)=(p^{\alpha}_2)^{-1}$ in the plane $C_2$ going
 in the opposite direction, replacing (+) by (-) and
the paths $\kappa,\kappa_1,\kappa_2$ by the similar paths $\delta,
\delta_1,\delta_2$: $$(p^1_2)^{-1}\sim
a^*_2\delta_1^{-1},(p^2_2)^{-1}b^*_2\sim\delta_2^{-1},
(p^0_2)^{-1}\sim a^*_2b^*_2\delta_2^{-1}$$. We have
$$\kappa_1=\delta_2,\kappa_2=\delta_1,\delta=\delta_1\delta_2^{-1}=\kappa^{-1}$$
We use the new basis $a^*_1,b^*_1,a^*_2,b^*_2$ defined above,
dropping the measure $m$ and signs $\pm$, as it was indicated
above.  {\bf All formulas are written in the new $m$-dependent basis}.

 Performing very simple
multiplication of paths, we obtain following formulas in the group
$\pi_1(V,s^+_1)$:

\begin{lem}
The homotopy types of all nonnegative almost transversal
two-street passes in the positive time direction starting and
ending in the segment $s^+_1=s^-_2$,  including the trajectory
passes, are equal to the following list of values of the elements
$\phi_{\alpha\beta}\in \pi_1(V,s^+_1)$ where $\alpha=1,0,2$,
$\beta=1',0',2'$ and $\alpha\beta$ means $\phi_{\alpha\beta}$:

 $$11'\sim a^*_1\kappa^{-1}(a^*_2)^{-1}; 10'\sim a^*_1b^*_1(a^*_2)^{-1};
  12'\sim b^*_1(a^*_2)^{-1}$$

$$01'\sim a^*_1(b^*_2)^{-1}(a^*_2)^{-1}; 00'\sim a^*_1b^*_1\kappa
 (b^*_2)^{-1}(a^*_2)^{-1}; 02'\sim b^*_1\kappa(b^*_2)^{-1}(a^*_2)^{-1}$$

$$21'\sim a^*_1(b^*_2)^{-1}; 20'\sim a^*_1b^*_1\kappa(b^*_2)^{-1}; 22'\sim
b^*_1\kappa(b^*_2)^{-1}$$
 For the negative time direction the two-street passes have homotopy classes equal to
$\phi_{\alpha\beta}^{-1}$ in the same group
$\pi_1(V,s^-_2)=\pi_1(V,s^+_1)$

 For the
Topological Types I--VI of Foliations following homotopy classes
of  almost transversal  two-street passes have nonzero measure:
$$\phi_{\alpha\beta}=p^{\alpha}_1xp^{\beta}_2,\phi^*_{\alpha\beta}=
(p^{\alpha}_1)^{-1}x(p^{\beta}_2)^{-1}$$  (here $x\geq 0$ means
positive transversal shift along the segment $s$):

 Type (I): All 9 classes $\phi_{\alpha\beta}$, and all  classes
 $\phi^*_{\alpha\beta}$ except $11',10',01',00'$

Type (II): All classes $\phi_{\alpha\beta}$ except $22'$, and all
classes $\phi^*_{\alpha\beta}$ except $11',10',01'$

Type (III): All classes $\phi_{\alpha\beta}$ except $22',02'$, and
all classes $\phi^*_{\alpha\beta}$ except $11',01'$

Type (IV): All classes $\phi_{\alpha\beta}$ except  $20',22'$, and
all classes $\phi^*_{\ alpha\beta}$ except $11',10'$

Type  (V): All classes $\phi_{\alpha\beta}$ except $22',20',02'$,
and all classes $\phi^*_{\alpha\beta}$ except $11'$

Type (VI): All classes $\phi_{\alpha\beta}$ except
$20',22',00',02'$, and all 9 classes $\phi^*_{\alpha\beta}$

\end{lem}

So we are coming to the following

{\bf Combinatorial Model of Foliation/Dynamical System}: {\it It involves\\
1.The semigroup $S_{\sigma,\tau}$ and its representation in the fundamental group. All
 positive finite words $$R=R_{j_1}R_{j_2}\ldots R_{j_N}\in \pi_1(V)$$
are written in the new $m$-dependent  Transversal Canonical Basis
$a^*_1,b^*_1,a^*_2,b^*_2$
where every symbol $R_j$ is equal to one of the elements
$\phi_{\alpha\beta}\in \pi_1(V),j=1,2,3,4,5$ found above. 2.The  permutation
$\sigma$ and measures of the sub-segments
$\tau_q,q=1,2,3,4,5, \sum_q\tau_q=m$. They  represent
 the types
$I,II,III,IV,V,VI$ described above. The permutation $\sigma$ defines a
broken isometry $i_{\sigma}:s\rightarrow s$ well-defined for the
inner points of the sub-segments $\tau_q$. The measure is extended to
the whole semigroup according to the rules formulated above.
3.Every finite word  $R\in S_{\sigma,\tau}$ defines a connected
strip of the nonseparatrix  trajectories of the Hamiltonian System
 corresponding to our data, with the transversal measure
$\tau_R$, starting and ending in the transversal segment $s^+_1$.
It defines a  closed transversal curve $\gamma_R$ with
transversal measure equal to the shift of the end $r(R)$. It is positive if $r(R)<0$,
 and negative
 if $r(R)>0$. The transversal measure of
these transversal curves are equal to $-r(R)$. The  infinite
trajectories are presented by the infinite sequences of the
symbols $R_j$ such that every finite piece of the sequence is
nonzero as an element of the semigroup $S_{\sigma,\tau}$. The
measure on the set of trajectories is defined by the transversal
measure of the finite words in $S_{\sigma,\tau}$.
The basic measurable sets are ''cylindrical'' (i.e. they consists
of all trajectories with the same finite word $R$, and measure is
equal  to $\tau_R$).}

Therefore the Step 2 is realized.

 {\bf Problem. How to classify simple closed transversal curves?}
 We consider this problem for the special case of torus with obstacle.
 Let a straight line flow on the 2-torus
with irrational rotation number  and fixed positive
transversal canonical basis $a,b^{-1}$ is given  as
above. Every indivisible homology class in $H_1(T^2,Z)$ can be
realized by the non-selfintersecting closed curve transversal to
foliation. We can see that {\bf the semigroup of positive closed
transversal curves has infinite number of generators containing
the elements with arbitrary small transversal measure}. The same
result remains true after removal of one point $T\rightarrow T^*=T^*_0$,
but the semigroup became non-abelian. Remove now  transversal
segment $s$ with positive measure $m>0$ from the torus
$T^*_0\rightarrow T^0_{(m)}$. As we shall see below, the situation
changes drastically: the semigroup became finitely generated with two free
 generators.

Start with Fundamental Domain in $C$
consisting of the union of the three streets $D_m=p^1\bigcup p^0\bigcup p^2$.
 We use an extension $D_m\subset D'\subset C$
of it adding two more streets-- one more copy of the street $p^1$ over
$p^2$ and one more copy of the street $p^2$ over $p^1$ (see Fig 16).
Every positive simple transversal curve can be
expressed in the $m$-dependent basis $a^*,b^*\in \pi_1(T^2_m)$. Its homology class is
$[\gamma]=ka^*+l(-b^*)\in H_1(T^*_m,Z), k>0,l>0, (k,l)=1$. We take
 the cycles $a=a^*,b^*=b^{-1}$ as  basic positive
closed transversal curves.

 Every non-selfintersecting positive closed transversal
curve $\gamma\subset T^*_m=T^2 minus (s)$ crosses fundamental domain and especially
the street $p^0$
 exactly $k+l$ times from the right to the
left. We denote its segments by $t_j,j=1,...,k+l$ ordering them
from the segment $s^+_1=s$ up. It crosses the
street $p^2$ exactly $k$ times, entering the domain $p^0$ in the
points $y_1,...,y_k$; they are naturally ordered by height. It
crosses  $p^1$ exactly $l$ times, entering the domain
$0$ from $p^1$  in the points $y_{k+1},...,y_{k+l}$,  We  denote the
 points on the equivalent segments from the left side of fundamental domain by the
same figures with symbol ': on the (lower) left
side of the street $p^0$  there are the points $y'_{k+1},...,y'_{k+l}$;
the points $y'_1,y'_2,...,y'_k$ we have on the
 segment separating the upper  copy of the street
$p^2$ from the street $p^0$. Assume that $k>l$.
The sequence of segments
ordered by height,is following:

$$I:2\rightarrow 0\rightarrow 1:t_1=[1,(k+1)'],...,t_l=[l,(k+l)']$$

 $$II:2\rightarrow 0\rightarrow 2:t_{l+1}=[l+1,1'],...,t_{k}=[k,(k-l)']$$

$$III:1\rightarrow 0\rightarrow 2:t_{k+1}=[k+1,(k-l+1)'],...,t_{k+l}=[k+l,k']$$

We assign to every segment following homotopy class depending on
the group I,II,III and their product to the whole curve $\gamma$:

 $$\phi:t_j\rightarrow b', t_j\in (I);
\phi:t_j\rightarrow a',t_j\in (II),(III)$$

$$\phi(\gamma)=\phi(t_{q_1})...\phi(t_{q_{k+l}})\in\pi_1(T^*_m)$$
where $\gamma\sim t_{q_1}...t_{q_{k+l}}$ in the natural order
along the curve. As a result of the previous lemmas, we obtain
following

\begin{thm}
The invariant $\phi (\gamma)$ describes  homotopy classes of all
closed non-selfintersecting positive transversal curves in $T^*_m$
as a positive words written in the free group with two
$m$-dependent generators $a'=a^*,b'=(b^*)^{-1}$ whose transversal
measures are smaller than $m$ . In particular, every indivisible
homology class $k[a']+l[b'],k>0,l>0,$ has  $m$-dependent simple
 positive transversal representative as
 a positive word unique up to cyclic permutation (i.e.
it defines exactly $k+l$ equivalent positive words). This word is calculated
above in the three street model. Every such curve can
be taken as a part of new Transversal Canonical Basis in $T^*_m$.
\end{thm}

{\bf Description of all Transversal Canonical Bases on the torus with obstacle}.
We can see that every integer-valued $2\times
2$-matrix $T,\det T=1,$ with positive entries $k,l,p,q\geq 0$
$$[A]=T(a')=k[a']+l[b'],[B]=T(b')=p[a']+q[b']$$
determines finite number of different   transversal canonical
bases $A,B\in \pi_1(T^*_m)$. They are  represented by the
curves $A,B$ crossing each other transversally in one point
and representing the homology classes $T(a'),T(b')$. Consider
 the  segments of both these curves in the street $p^0$
as above. Let the  curves $A,B$ are represented by the ordered
sequences of pieces $t_1,...t_{k+l}$ and $t''_1,...t''_{p+q}$
correspondingly going from the right to the left side. We require
existence of  one intersection point for the selected pair
$t_i\bigcap t''_j\neq \emptyset$. All other intersections should be
empty. Every such configuration determines transversal canonical
basis $A,B$. Its equivalence class  is completely determined
by the relative order of segments $t,t''$ taking into account that
$t_i$ crosses $t''_j$ only once and $t_i$ is ''higher'' than
$t''_j$ from the left side). Starting from the selected point
$t_i\bigcap t''_j$, we apply the procedure described  in the
theorem. It gives us two positive words $A,B$ in the free group
$F_2$. By definition, the map $\hat{T}:a'\rightarrow
A,b'\rightarrow B$ defines a lift from homology to fundamental
group, for every pair of words $A,B$ constructed in that way.
Consider first a {\bf Reducible}  case where the left ends of the
segments $t_i$ and $t''_j$ crossing each other are located on the
same connected part of boundary (i.e. where the street $p^0$ meets the same street
$p^1$ or $p^2$). If it is the street $p^1$,  than the both words
$A$ and $B$ start with the same letter $b'$. We can deform our
crossing point along the cycle $b'$. After this step we are coming
to the conjugated pair $A',B'$ where $b'$ is sent to the end:
$$(A=b'\tilde{A},B=b'\tilde{B})\rightarrow (A'=\tilde{A}b',B'=\tilde{B}b')$$.
 If the pair $t_i,t''_j$ ends up in the
 street $p^2$, we replace $b'$ by
$a'$. After the series of such steps, we  are coming to the case
where $t_i$ ends in $p^1$,  $t''_j$ ends in $p^2$. This process cannot
be infinite because the words $A,B$ are not  powers of the same
word. So this process ends. We call  {\bf reduced} the case
where the process ended up.

 {\bf Example}: The pair of words
$A=b'a'b'a'b', B=b'a'$ requires five conjugations
to arrive to the reduced case.

In the final reduced state
 the relation:
$ABA^{-1}B^{-1}=a'b'a'^{-1}b'^{-1}$
can be easily seen on the plane (with periodic set of segments
removed), looking on the 3-street decomposition of the plane.

 {\bf The Semigroup of unimodular $2\times 2$-matrices with
nonnegative integer-valued entries is free, with two generators
$T_1,T_2$} such that
$$T_1(a')=a'+b',T_1(b')=b',T_2(a')=a',T_2(b')=a'+b'$$
Their lifts to the automorphisms of the free group are following:
$$\hat{T}_1(a')=a'b',\hat{T}(b')=b';\hat{T}_2(a')=a',\hat{T}_2(b')=b'a'$$
They both preserve the word $a'b'a'^{-1}b'^{-1}$, so
all semigroup of nonnegative unimodular matrices preserves this
word. It is isomorphic to the semigroup of all positive
automorphisms of free group $F_2$ preserving this word.\\
 The semigroup $G$ of all positive automorphisms
of the free group $F_2$ in the alphabet $a',b'$ preserving the
conjugacy class of the word $\kappa=a'b'a'^{-1}b'^{-1}$, contains
following parts:
1.A free semigroup $G_{\kappa}$  of  transformations
preserving this word $\kappa$ exactly; It is isomorphic to the semigroup of
matrices $T$ with $\det T=1$ and nonnegative integer entries;
2. Every element    $T\in G_{\kappa}$ defines  finite number of
positive transformations  $T'\in G$ such that the corresponding
pair of positive words $A'=T'(a'),B'=T'(b')$ is simultaneously
conjugate to the words $A=T(a),B=T(b)$ for $T\in G_{\kappa}$. All these
conjugations are the simultaneous cyclic
permutation of the words $A,B$ removing the same letter from the left
end  and sending it to the beginning (until both words ends up with
the same letter). {\bf Total number of these conjugations is
equal to the sum of matrix elements of $T$ minus 2}. The natural
projection of semigroups $h:G\rightarrow G_{\kappa}$ is such that the
inverse image $h^{-1}(T)$ of each  matrix  $T$ contains
exactly $k+l+p+q-2$ positive automorphisms. In the example above
we have $k+l+p+q-2=5$. We clarified this question with I.Dynnikov.
I don't know where it is written.\\
{\bf The ''$m$-cutted Euclidean Algorithm''}: How to find $a',b'$ starting from
the arbitrary TCB $A,B$? Let $A,B\subset T^*_m, |A|+|B|>m$ and
 $|A|>|B|$. Take {\bf the minimal integer} $l_1\geq 0$ such that $$|A|-l_1|B|=|A_1|=r_1|B|
 <\max(m,|B|)$$
For the new pair $(|A_1|,|B_1|=|B|>|A_1|)$ take minimal integer $l_2\geq 0$ such that
$r_2|B|=|B_2|=|B_1|-l_2|A_1|<\max(m,|A_1|)$. Construct new pair  $(|A_2|=|A_1|, |B_2|)$
 and so on.
We have  sequence $l_1,l_2,\ldots$. The process stops at the first moment $q$ when
$|A_q|+|B_q|>m$ but  $\max(|A_q|,|B_q|)<m$. We put $|a'|=|A_q|,|b'|=|B_q|$. The transformation
$T=T_1^{l_1}T_2^{l_2}T_1^{l_3}T_2^{l_4}\ldots$ maps  basis $a',b'$ into $A,B$.
The Continued Fraction $|A|/|B|=l'_1+1/(l'_2+1/(l'_3+1/(...)))$ has exactly
 the same coefficients except the last one: Let
 $r_j=x_1\ldots x_j$ where $x_j=1/(l'_{j+1}+1/(l'_{j+2}+...))$ and $l'_{j+1}=[1/x_j]\in Z_+$.
 Until $r_{j-1}=x_1...x_{j-1}>m/|B|$, we have $l'_j=l_j$. Let $j$ be  the first number with property
  $r_j=x_1...x_j<m/|B|$.
 We define $l_{j+1}$ as a minimal integer such that $l_{j+1}>x_j^{-1}(1-m/(|B|x_1...x_{j-1}))$.
 At the same time, we have $l_{j+1}>x_j^{-1}$, so $l'_{j+1}\geq l_{j+1}$ and $j+1=q$.

\begin{center}
\mbox{\epsfxsize=5cm \epsffile{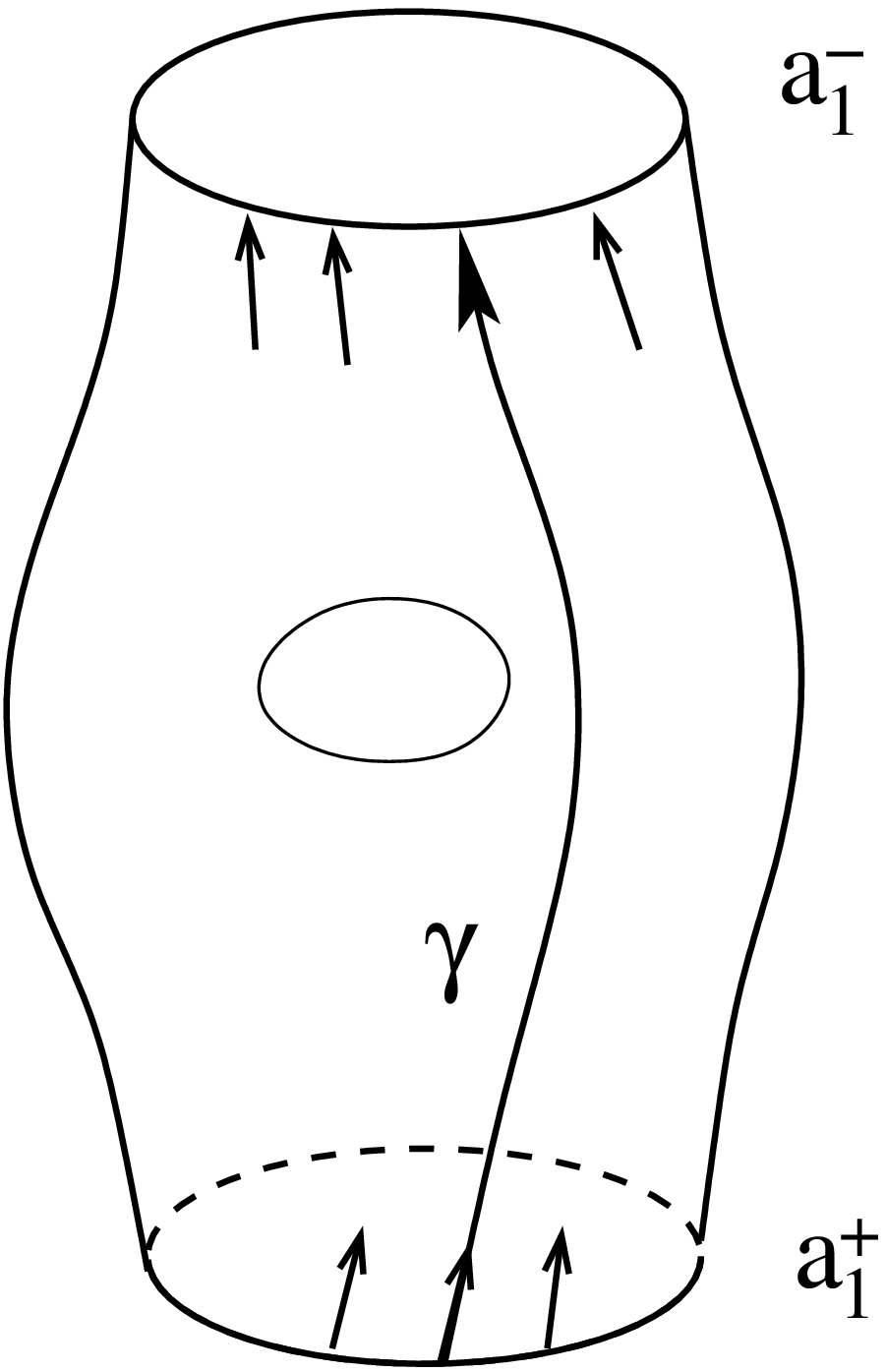}}

Fig 1: Constructing TCB: The first pair of cycles
\end{center}

\begin{center}
\mbox{\epsfxsize=4cm \epsffile{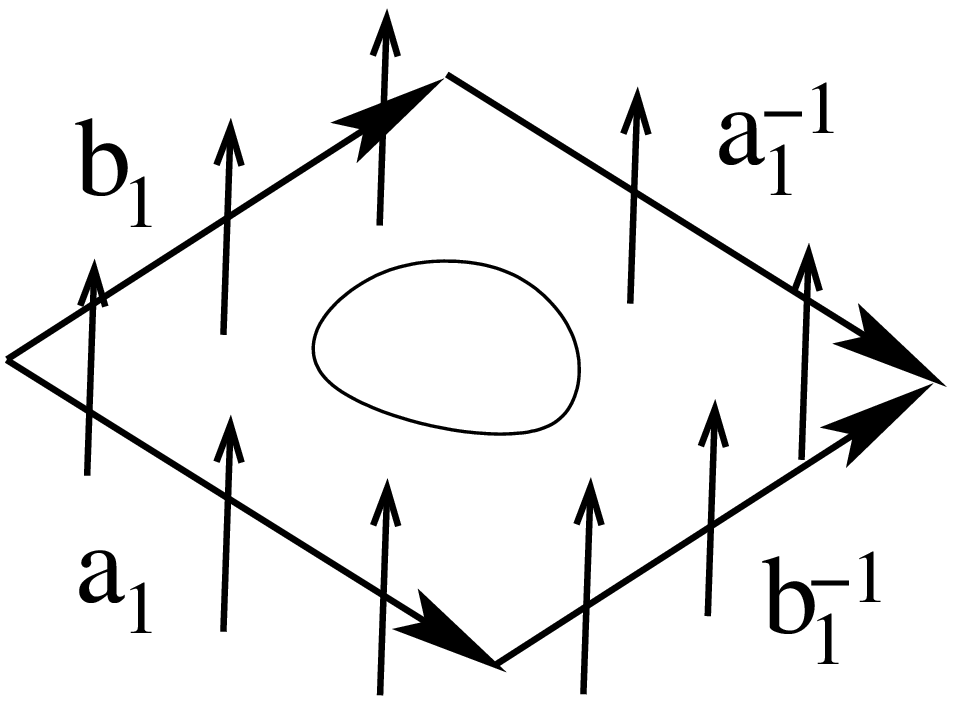}}

Fig 2: Cutting Riemann Surface along TCB.
\end{center}

\begin{center}
\mbox{\epsfxsize=10cm \epsffile{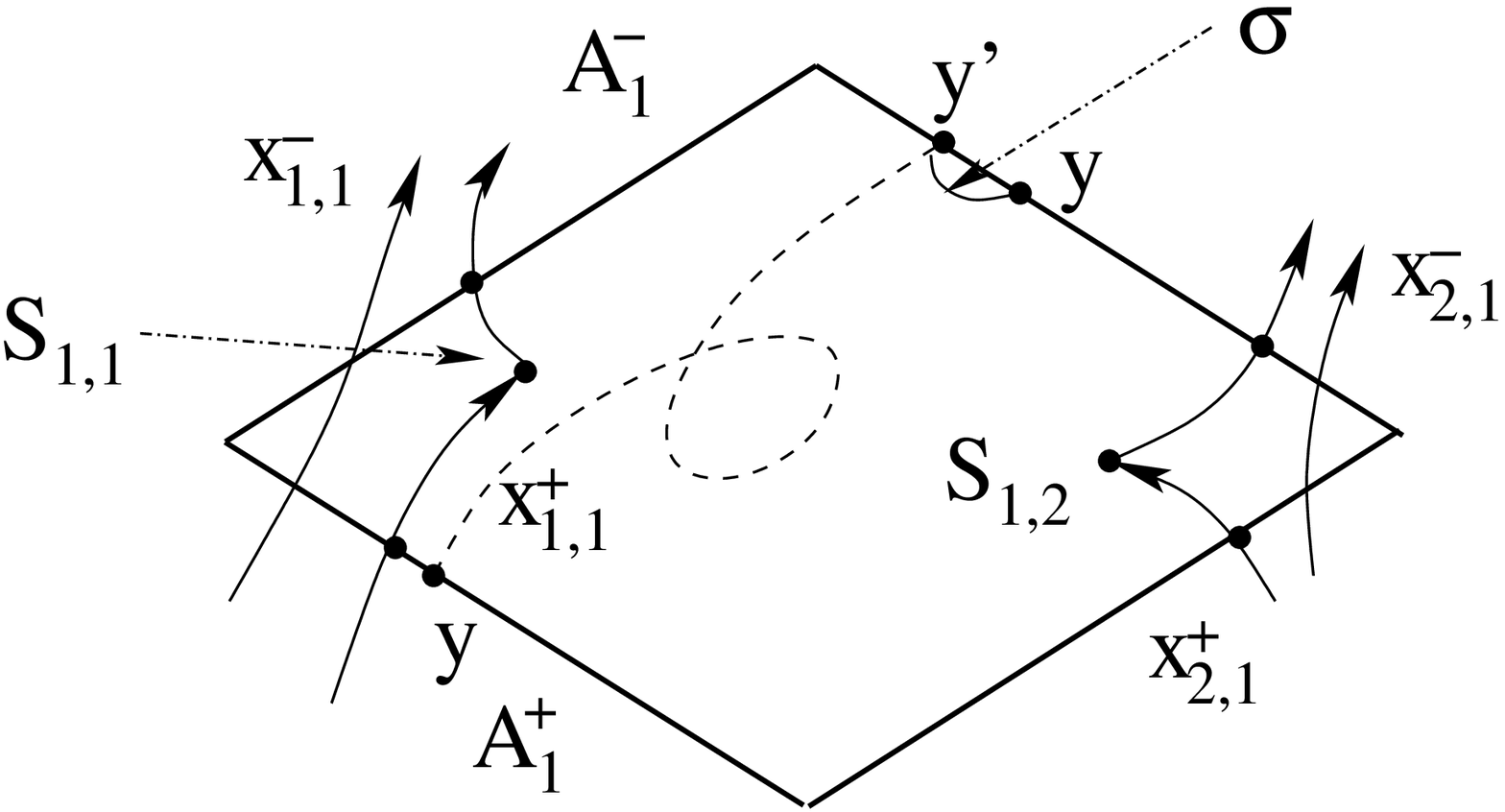}}

Fig 3: Boundary Saddles
\end{center}

\begin{center}
\mbox{\epsfxsize=10cm \epsffile{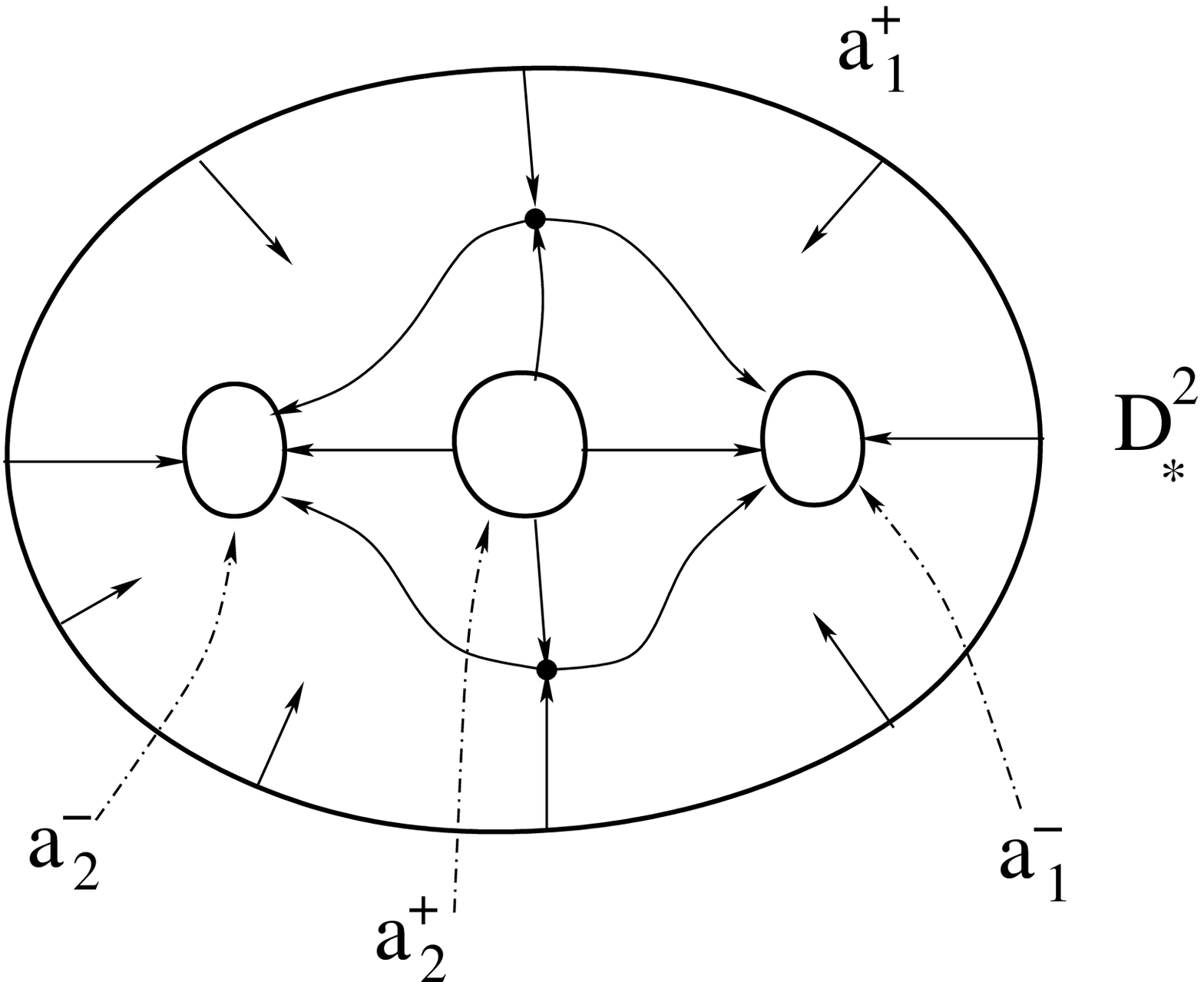}}

Fig 4a: Trajectory connections: $a_1^+\rightarrow a_1^-$,
         $a_2^+\rightarrow a_2^-$.

\end{center}

\begin{center}
\mbox{\epsfxsize=10cm \epsffile{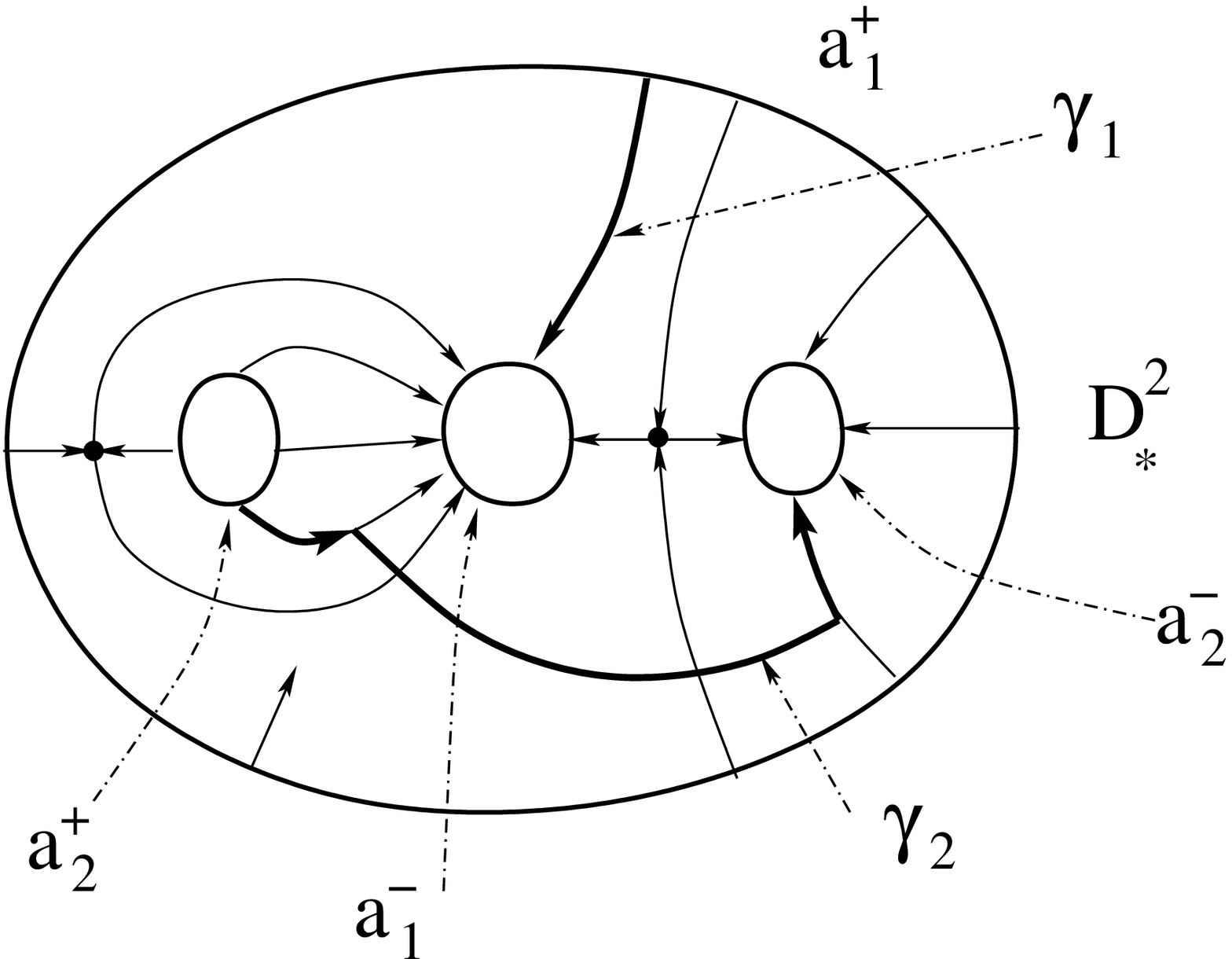}}

Fig 4b: Trajectory connection: $a_1^+\rightarrow a_1^-$.

Almost transversal curve $\gamma$ (boldface): $a_2^+\rightarrow
a_2^-$.

\end{center}

\begin{center}
\mbox{\epsfxsize=8cm \epsffile{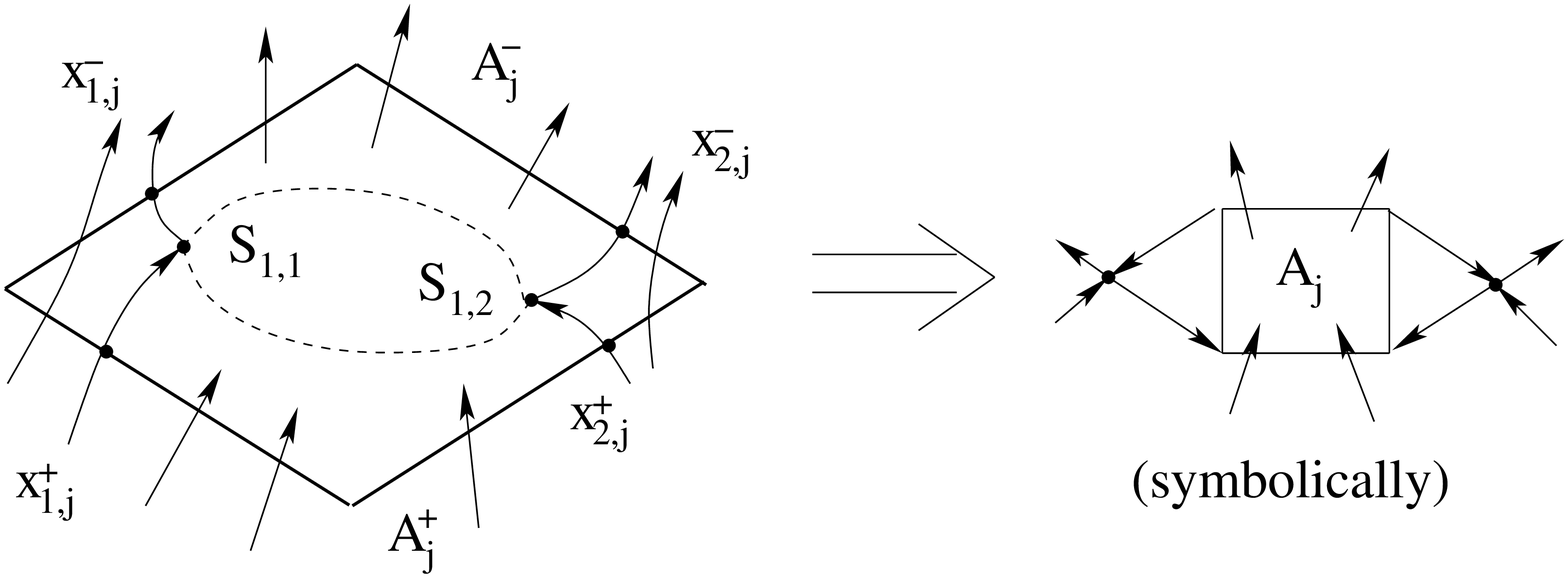}}

Fig 5: Transversal Cuts in the Riemann Surface
\end{center}

\begin{center}
\mbox{\epsfxsize=7cm \epsffile{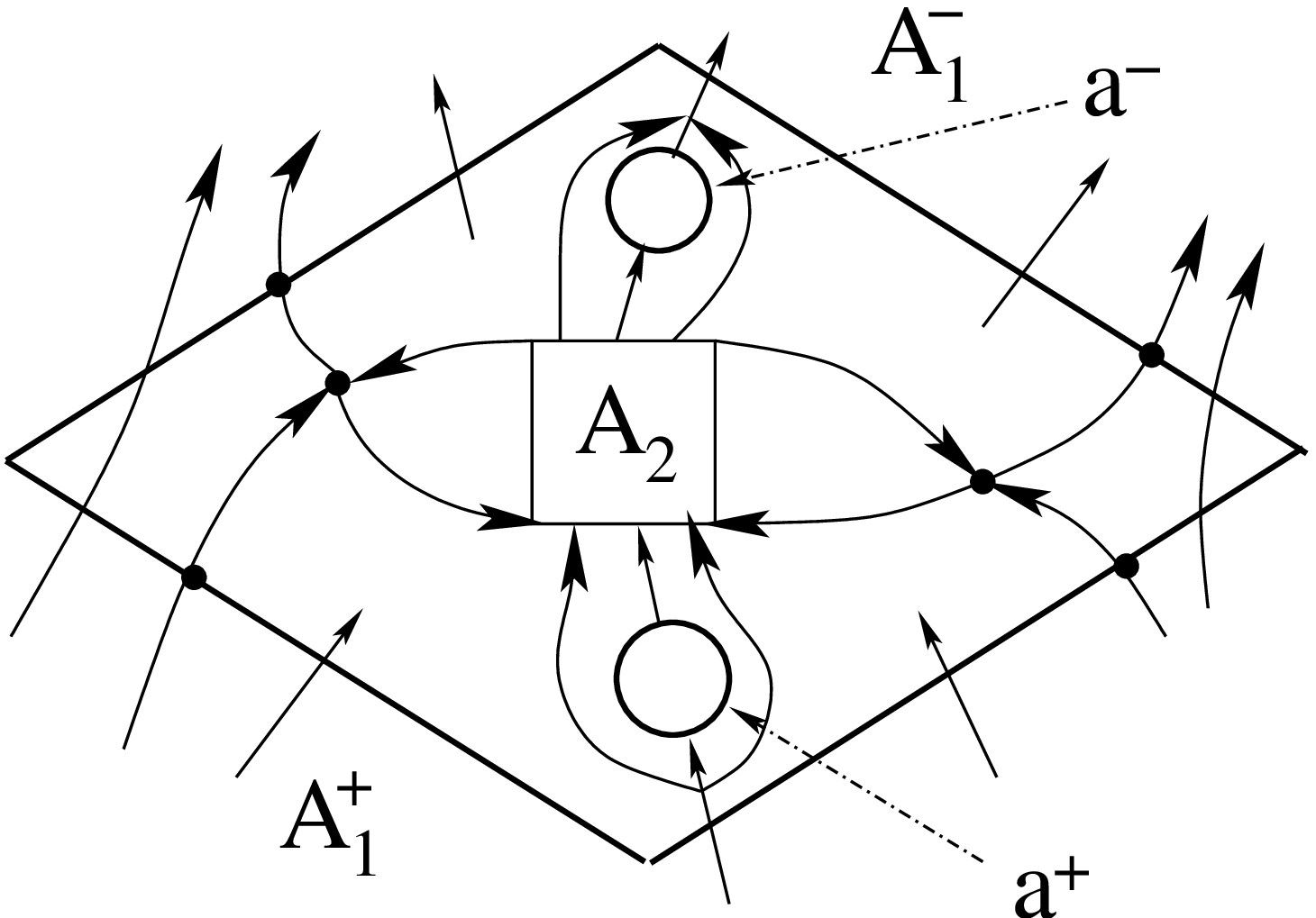}}

Fig 6: Non-extendable diagram of the type $T^2$.

\end{center}

\begin{center}
\mbox{\epsfxsize=9cm \epsffile{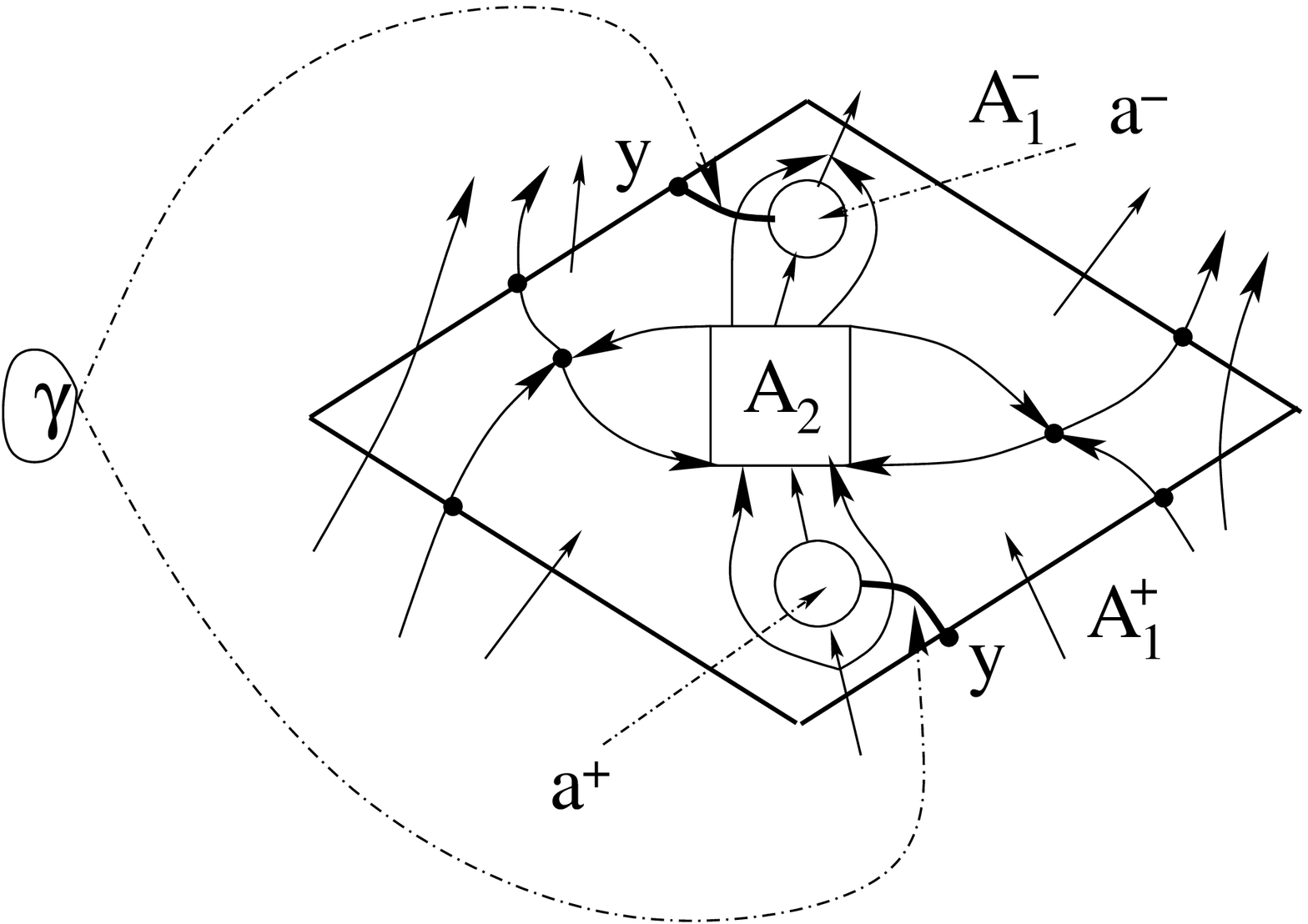}}

Fig 7: Reconstruction: new transversal curve $\gamma:$ $a_1^+\rightarrow
a_1^-$ is passing through the boundary.

\end{center}

\begin{center}
\mbox{\epsfxsize=6cm \epsffile{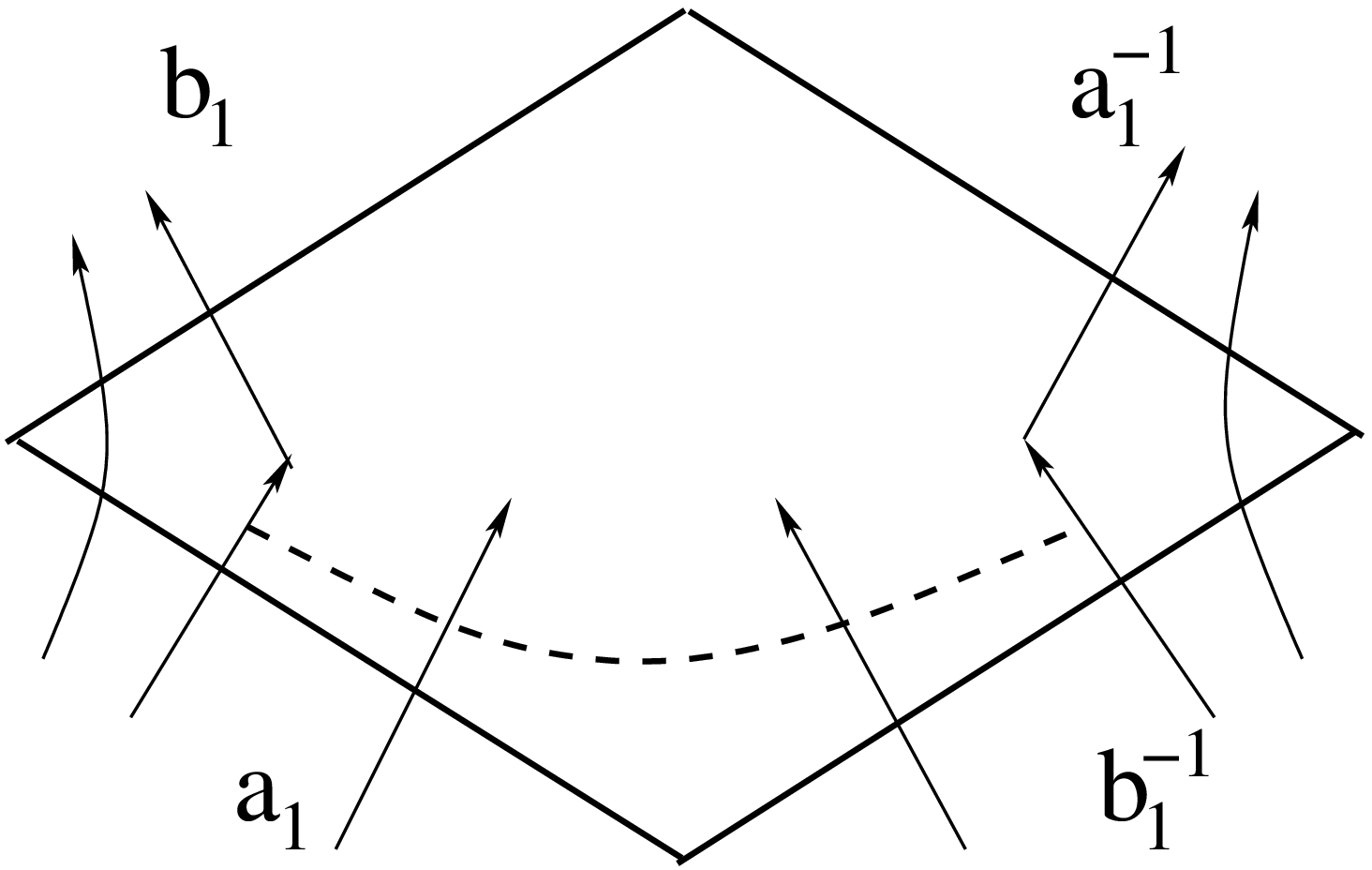}}

Fig 8: Transversal segment joining two saddles
\end{center}

\begin{center}
\mbox{\epsfxsize=6cm \epsffile{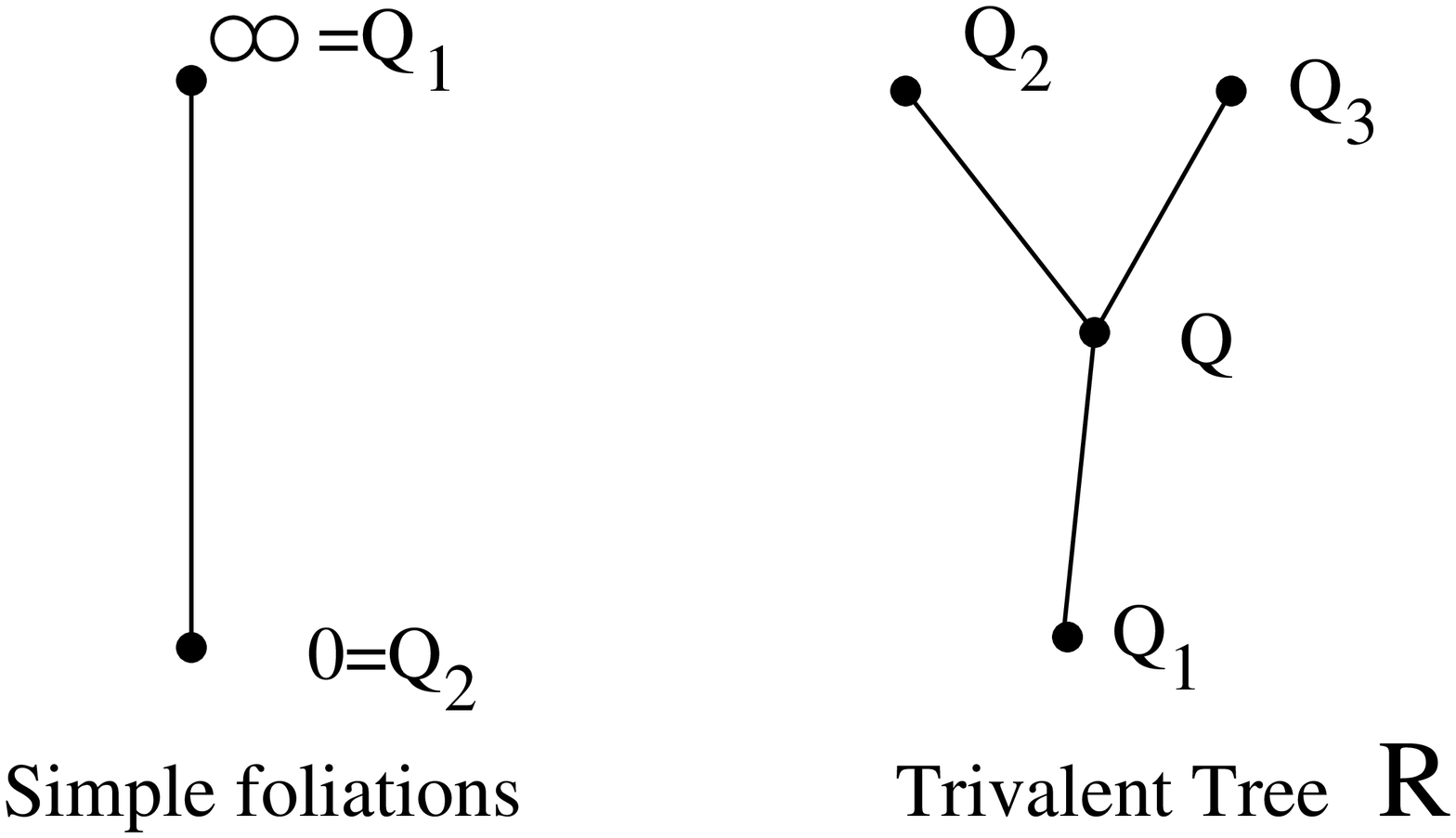}}

Fig 9c: Coding  Morse Functions on $S^2$ by Graphs
\end{center}

\begin{center}
\mbox{\epsfxsize=9cm \epsffile{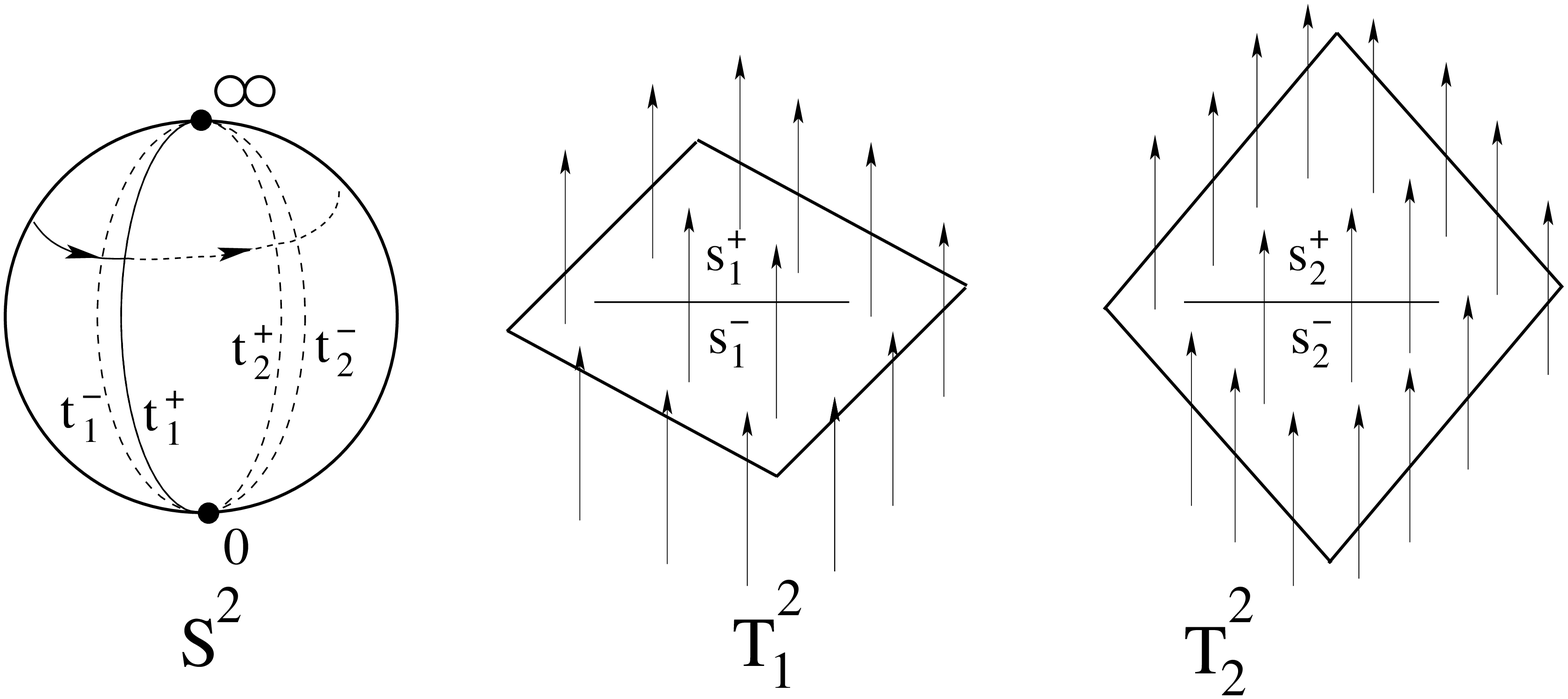}}

Fig 9a: The Building Data: $g=2$

\end{center}

\begin{center}
\mbox{\epsfxsize=9cm \epsffile{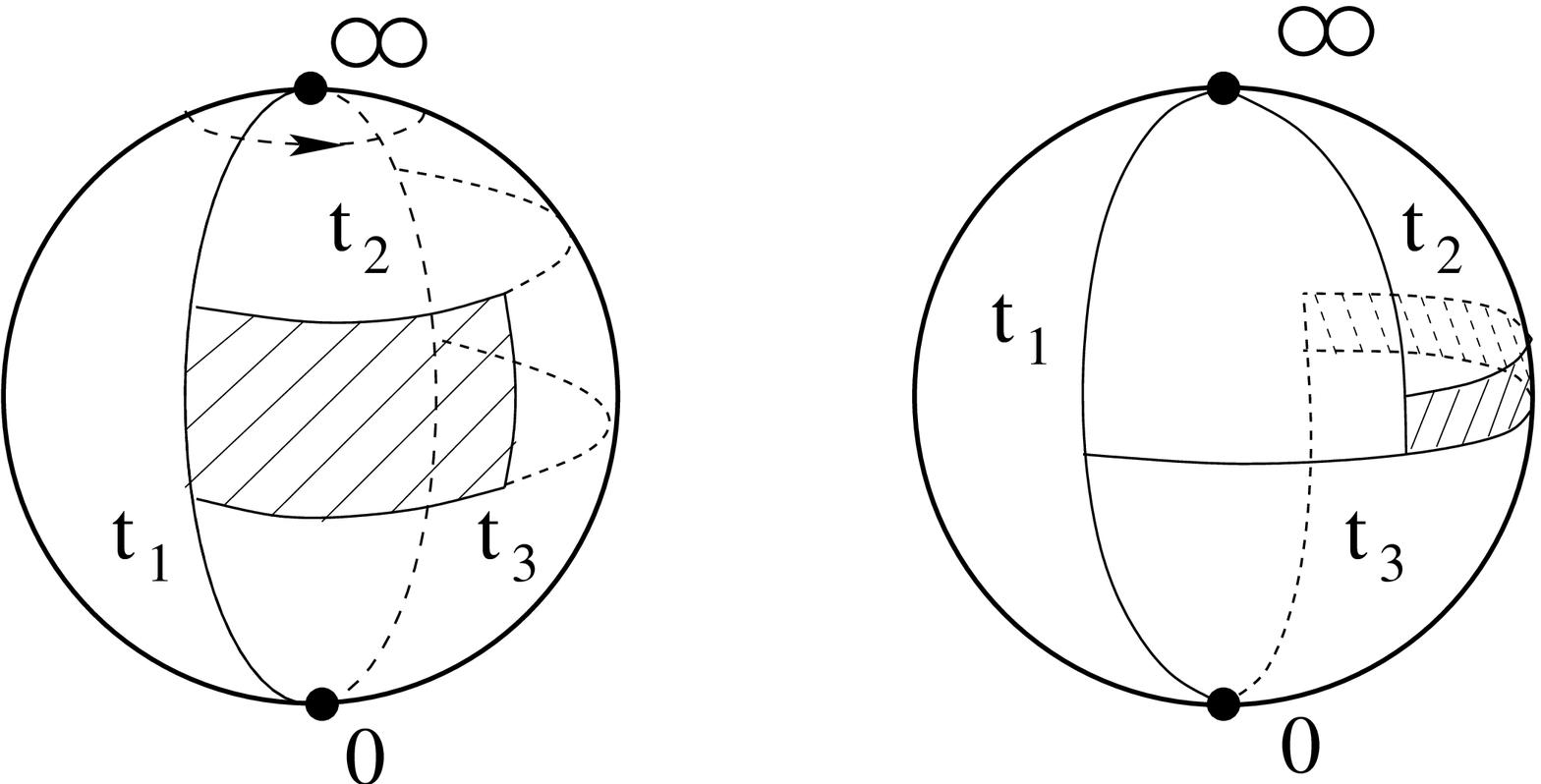}}

Fig 9b: The Corridors in $S^2$ ($g=3$)

\end{center}

\begin{center}
\mbox{\epsfxsize=6cm \epsffile{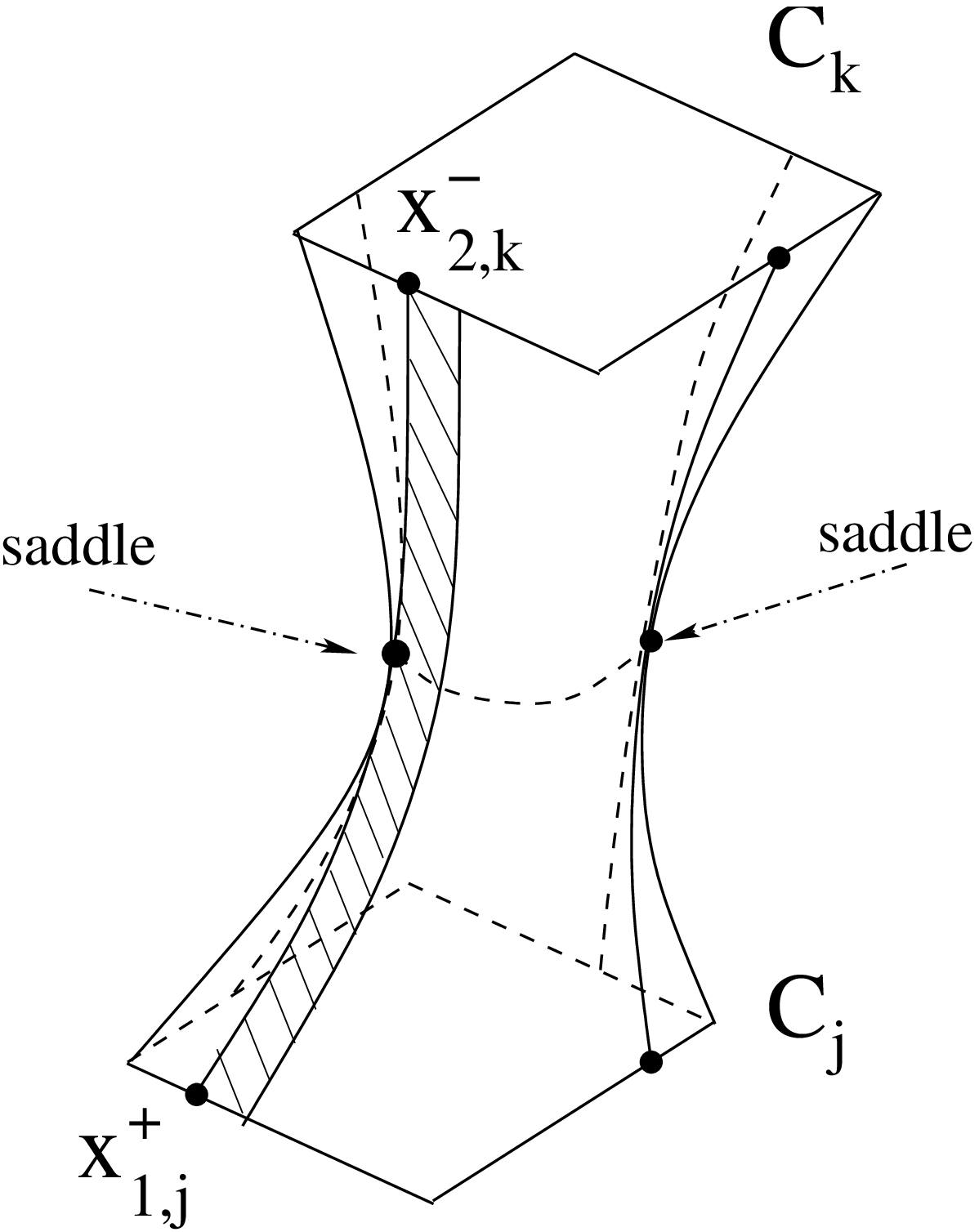}}

Fig 10: The Abelian Fundamental Domain ($g=2$)
\end{center}

\begin{center}
\mbox{\epsfxsize=12cm \epsffile{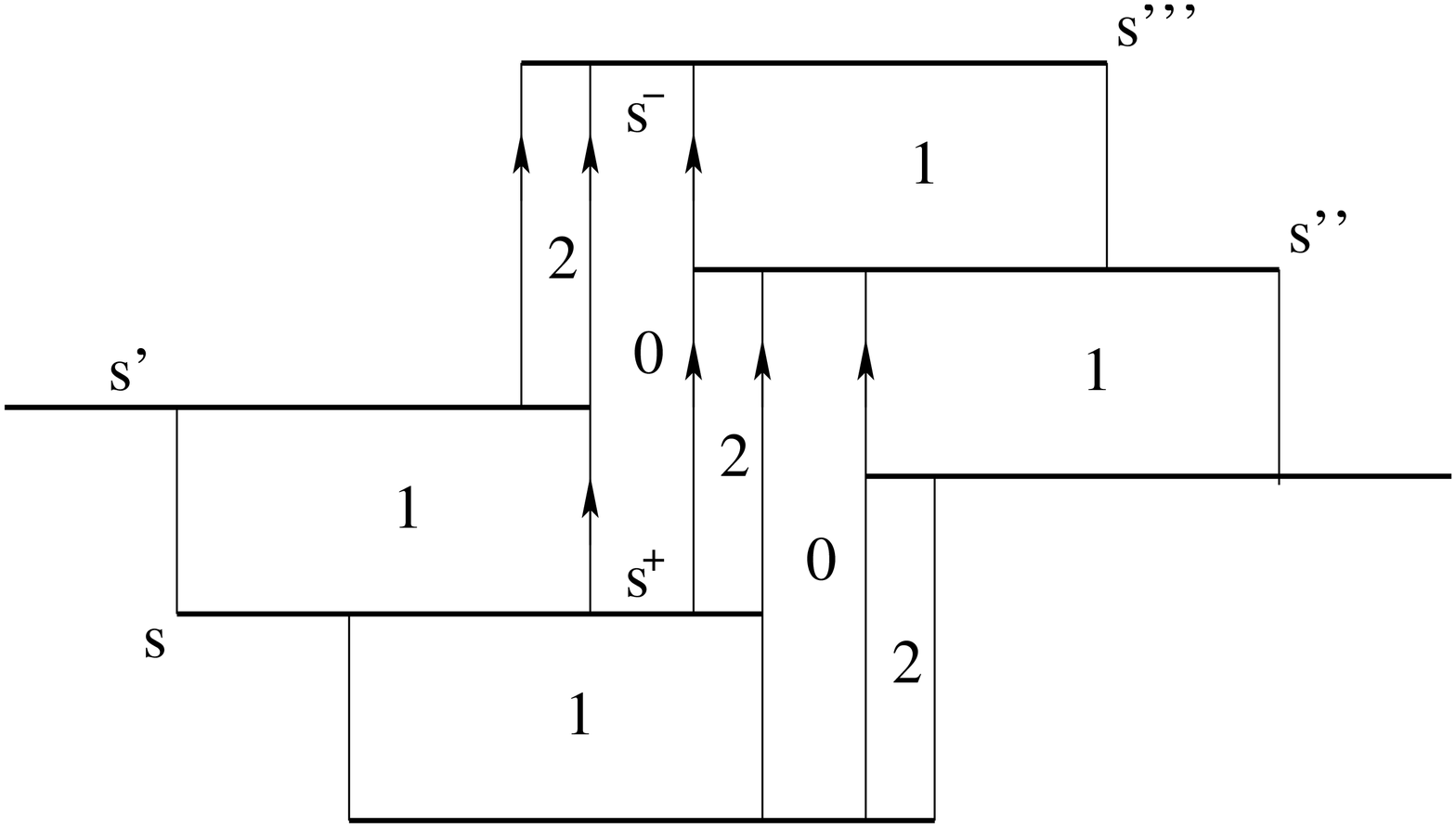}}

Fig 11: The 3-street model.

\end{center}

\begin{center}
\mbox{\epsfxsize=5cm \epsffile{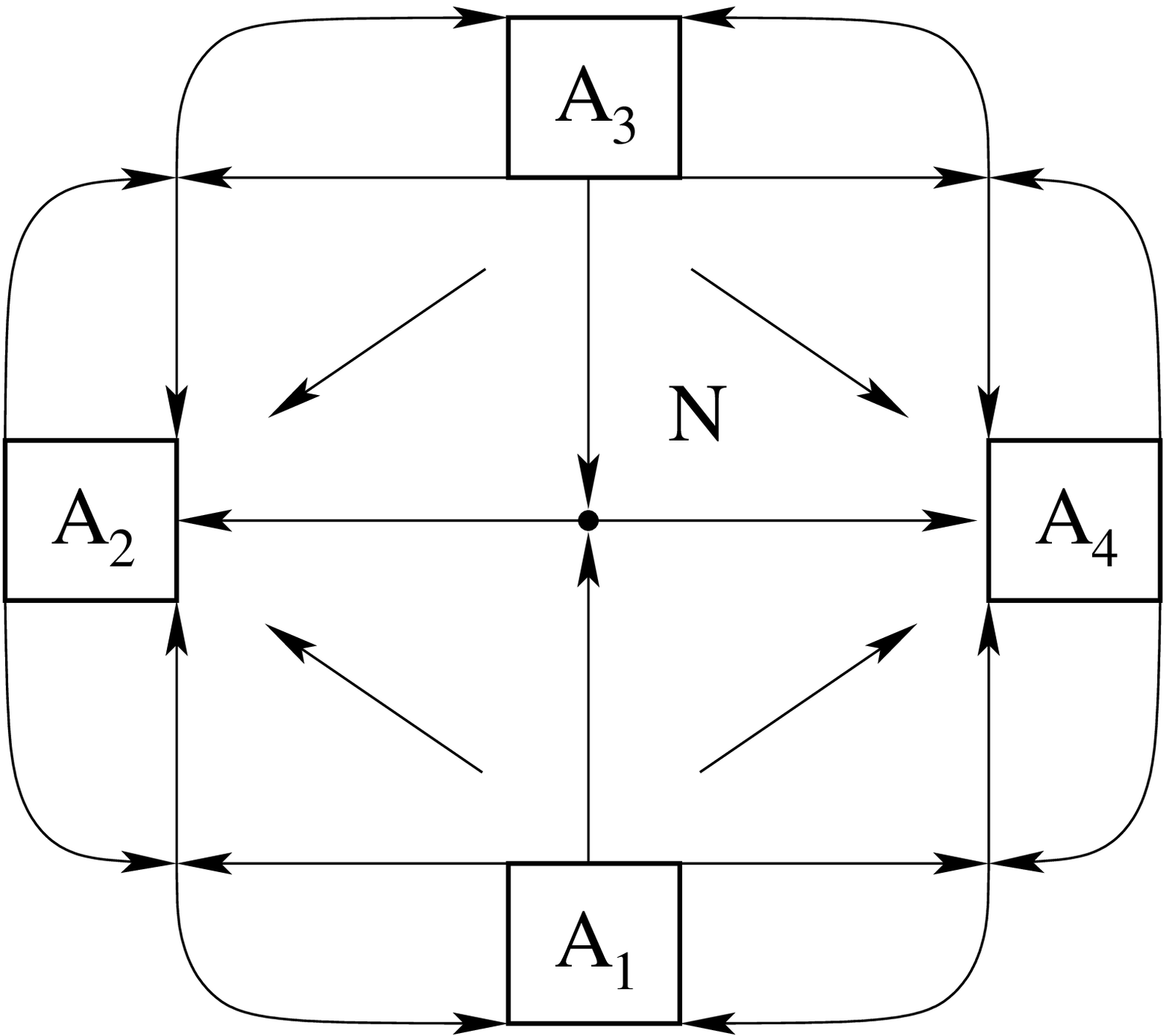}}\hspace{1cm}
\mbox{\epsfxsize=5cm \epsffile{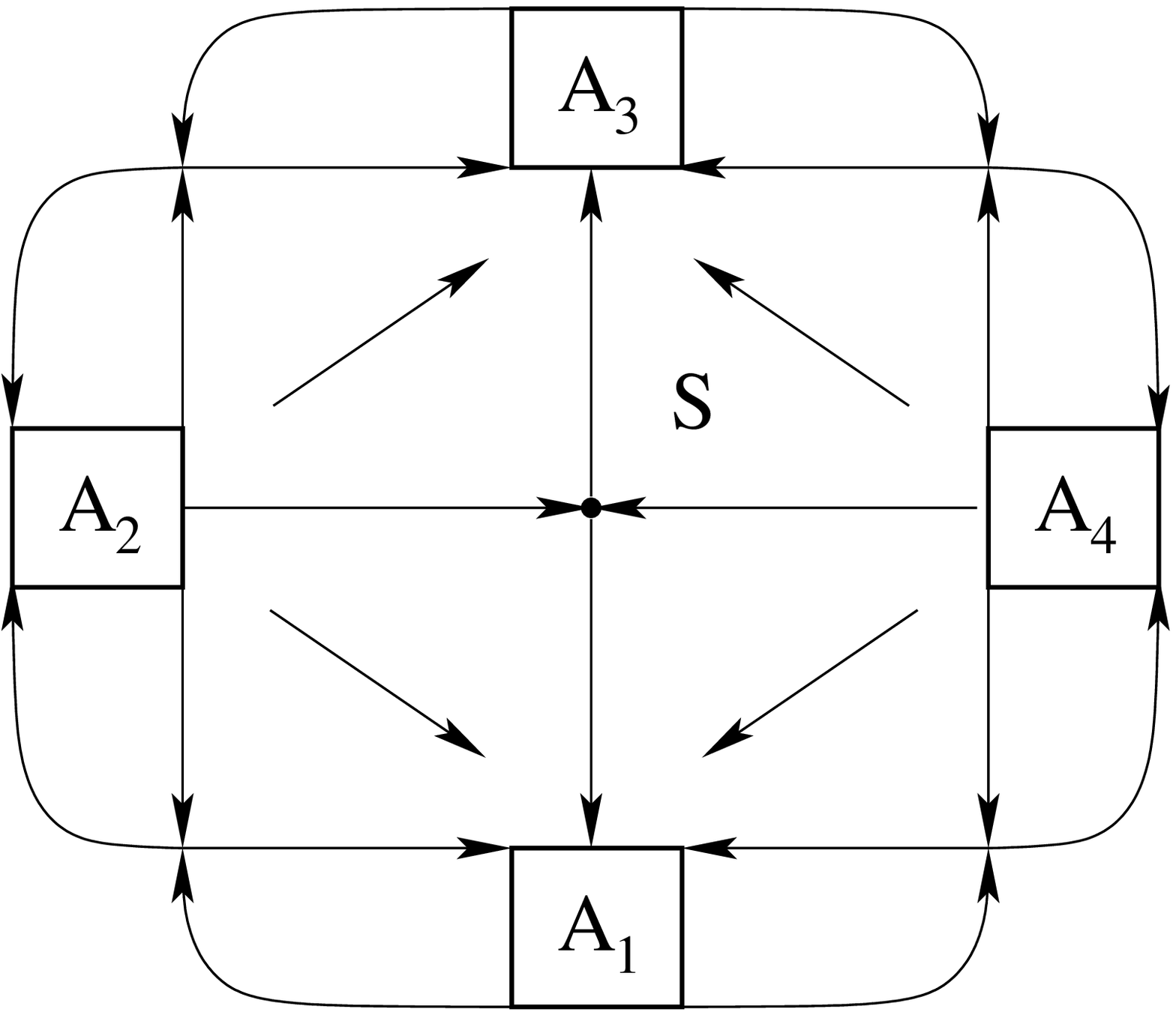}}

Fig 12: The Maximal Foliation:  North Hemisphere(left) and South Hemisphere(right)

\end{center}

\begin{center}
\mbox{\epsfxsize=6cm \epsffile{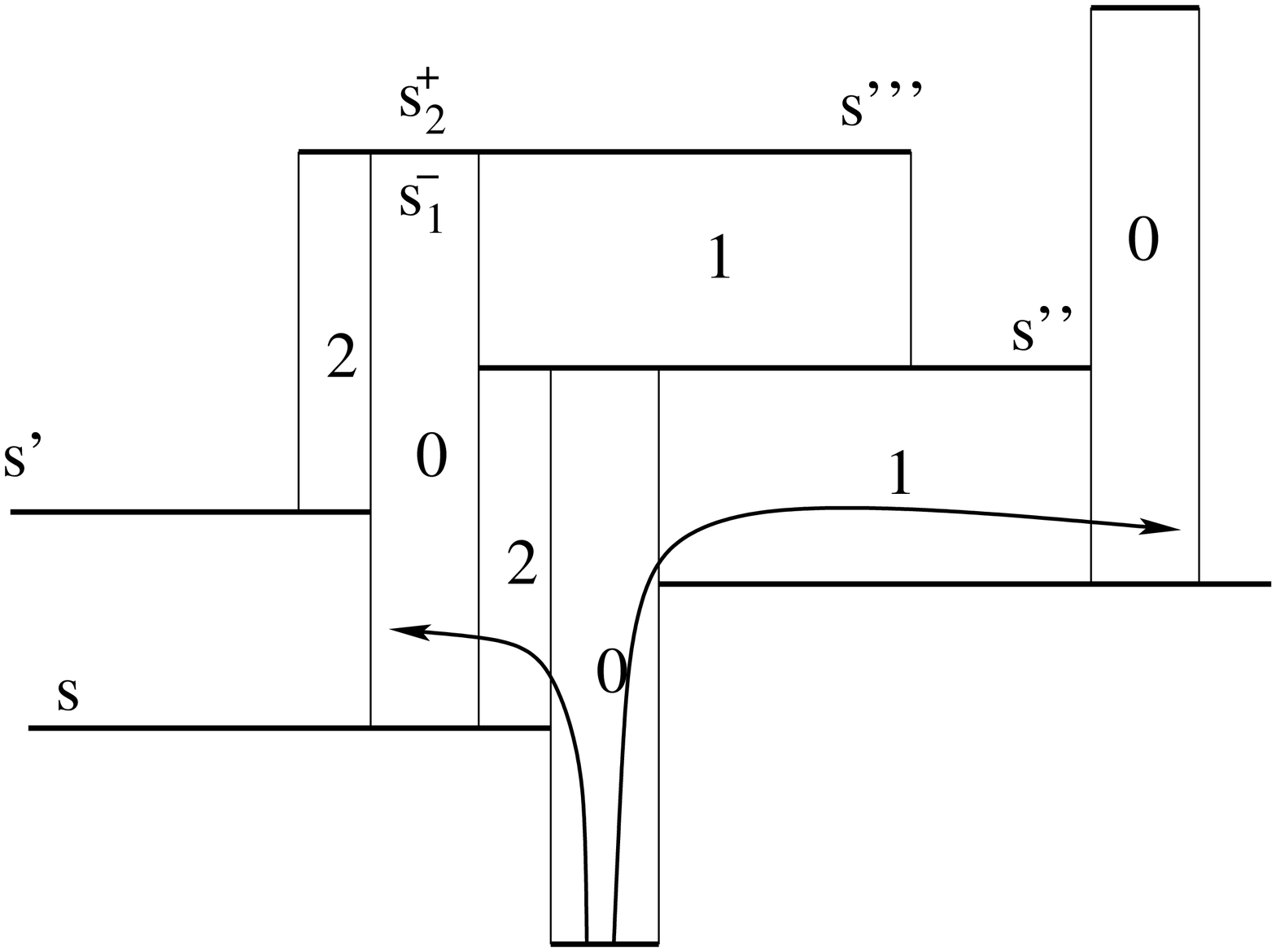}}

Fig 13: $a^*:0\rightarrow 2\rightarrow 0$, $b^*:0\rightarrow
1\rightarrow 0$.

\end{center}

\pagebreak

\begin{center}
\mbox{\epsfxsize=9cm \epsffile{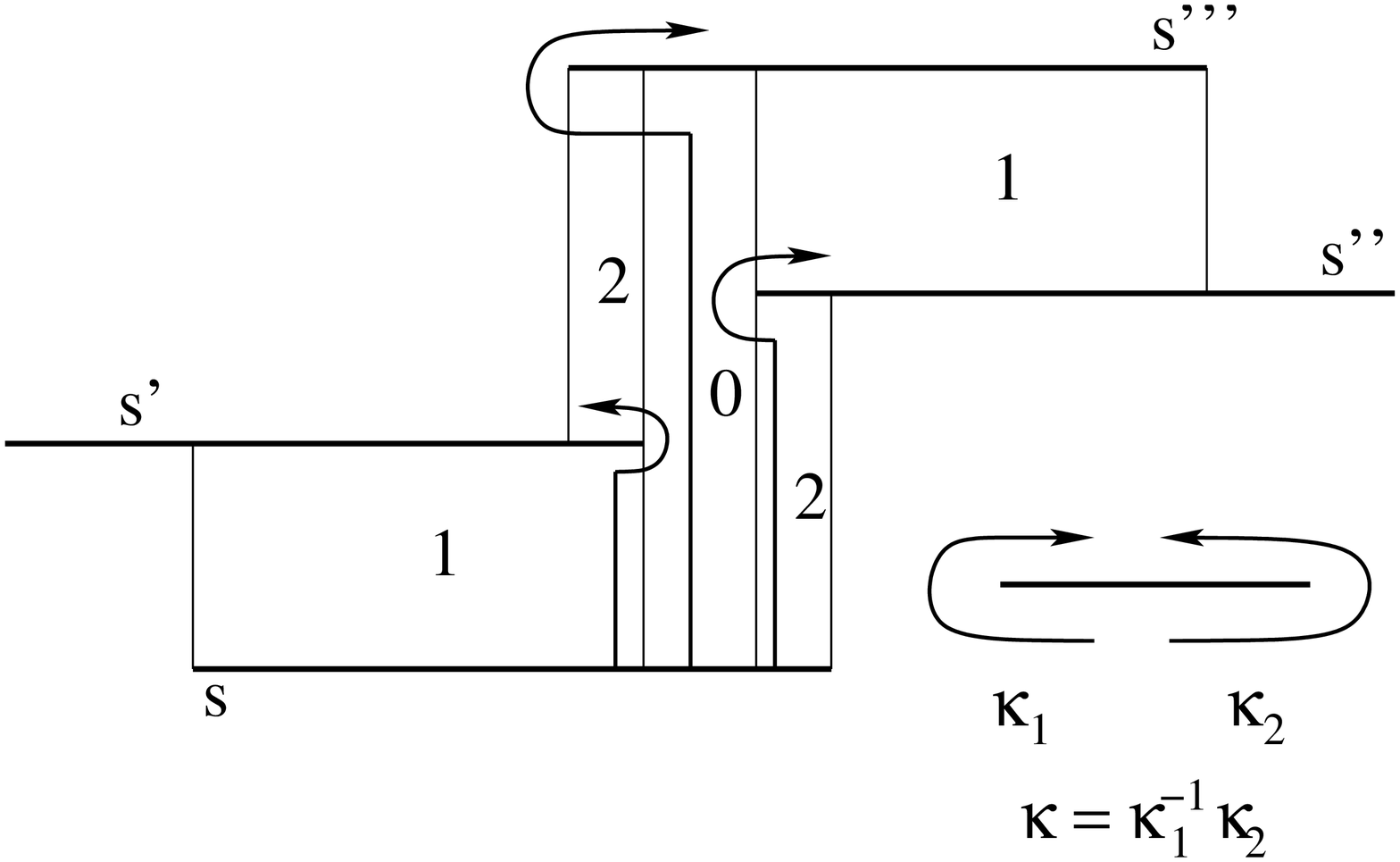}}

Fig 14. The elements of $\pi_1$ assigned to the streets
\end{center}

\begin{center}
\mbox{\epsfxsize=10cm \epsffile{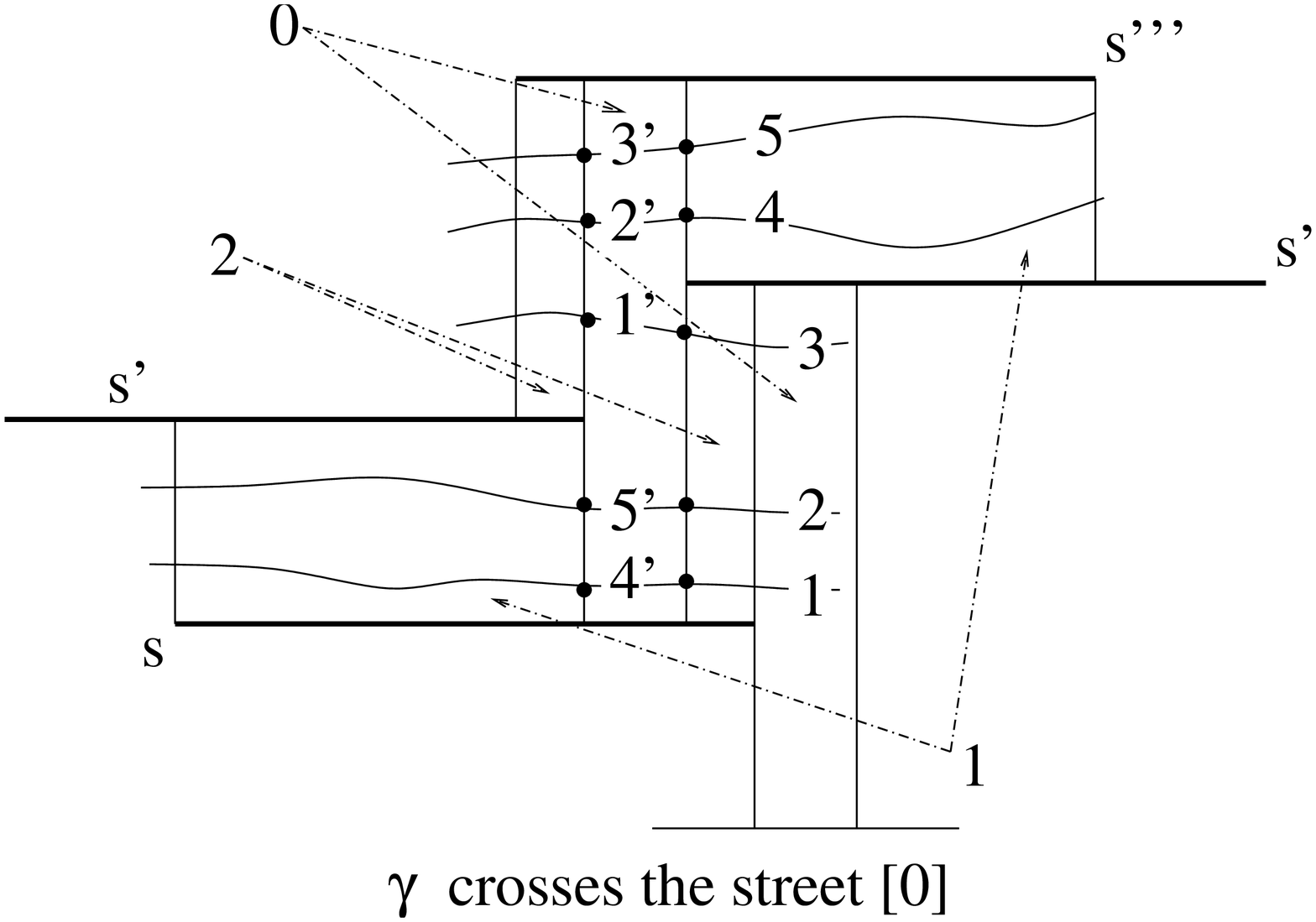}}

Fig 15:The simple transversal curve: $k=3$, $l=2$.

\end{center}

\end{document}